\documentstyle[11pt,psfig]{article}
 \newtheorem{theorem}{Theorem}[section]
 \newtheorem{remark}[theorem]{Remark}
 \newtheorem{lemma}[theorem]{Lemma}

\newtheorem{definition}[theorem]{Definition}
 \newenvironment{proof}{{\it Proof:\/}}{$\Box$\vskip 0.08in}
\newcommand{\mod}{{\mbox{ mod }}}
\newtheorem{problem}[theorem]{Problem}
\newtheorem{example}[theorem]{Example}
\newtheorem{corollary}[theorem]{Corollary}
\newtheorem{conjecture}[theorem]{Conjecture}

\newtheorem{proposition}[theorem]{Proposition}

\newcommand{\pct}[1]{}
\newcommand{\Psfig}[1]{{\mbox{$\ \ $}}}

\begin{document}
\vspace*{0.2in} 
\begin{quote} 
{\em 
Symmetry is a vast subject, significant in art and nature.
Mathematics lies at its roots, and it would be hard to find
a better one on which to demonstrate the working of the mathematical
intellect.
} 
{\footnotesize \ \ \ 
[Hermann Weyl, Symmetry \cite{We}]} 
\end{quote} 
\vspace{0.4in} 
 \begin{center} 
 {\LARGE\bf Symmetric knots and billiard knots. }
\end{center}
\begin{center} {\bf by J\'ozef H.Przytycki} \end{center}
\ \\

Symmetry of geometrical figures is reflected in regularities of
their algebraic invariants. Algebraic regularities are often
preserved when the geometrical figure is topologically deformed.
The most natural, intuitively simple but mathematically complicated,
topological objects are {\bf \it knots}.

We present in this papers several examples, both old and new, of
regularity of algebraic invariants of knots. Our main invariants are
the Jones polynomial (1984) and its generalizations.

In the first section, we discuss the concept of a symmetric knot,
and give one important example -- a torus knot. In the second section,
we give review of the Jones type invariants. In the third section, we gently
and precisely develop the periodicity criteria from the Kauffman bracket
(ingenious version of the Jones polynomial). In the fourth section,
 we extend the criteria to skein (Homflypt) and Kauffman polynomials.
In the fifth section we describe $r^q$ periodicity criteria using
Vassiliev-Gusarov invariants. We also show how the skein method may be 
used for $r^q$ periodicity criteria for  the classical (1928)
Alexander polynomial. In the sixth section, we introduce the notion of 
Lissajous and billiard knots and show how symmetry principles can be 
applied to these geometric knots. Finally, in the seventh section,
we show how symmetry can be used to gain nontrivial information about
knots in other 3-manifolds, and how symmetry of 3-manifolds is reflected
in manifold invariants.

\section{Symmetric knots.}
We analyze, in this  paper, symmetric knots and links, that is, links
which are invariant under a finite group action on $S^3$ or, more
generally, a 3-manifold. 

For example, a torus link of type $(p,q)$ (we call it $T_{p,q}$)
is preserved by an action of a group $Z_p\oplus Z_q$ on $S^3$ 
(c.f.~Fig.~1.1 and Example 1.1).

\ \\
\centerline{\psfig{figure=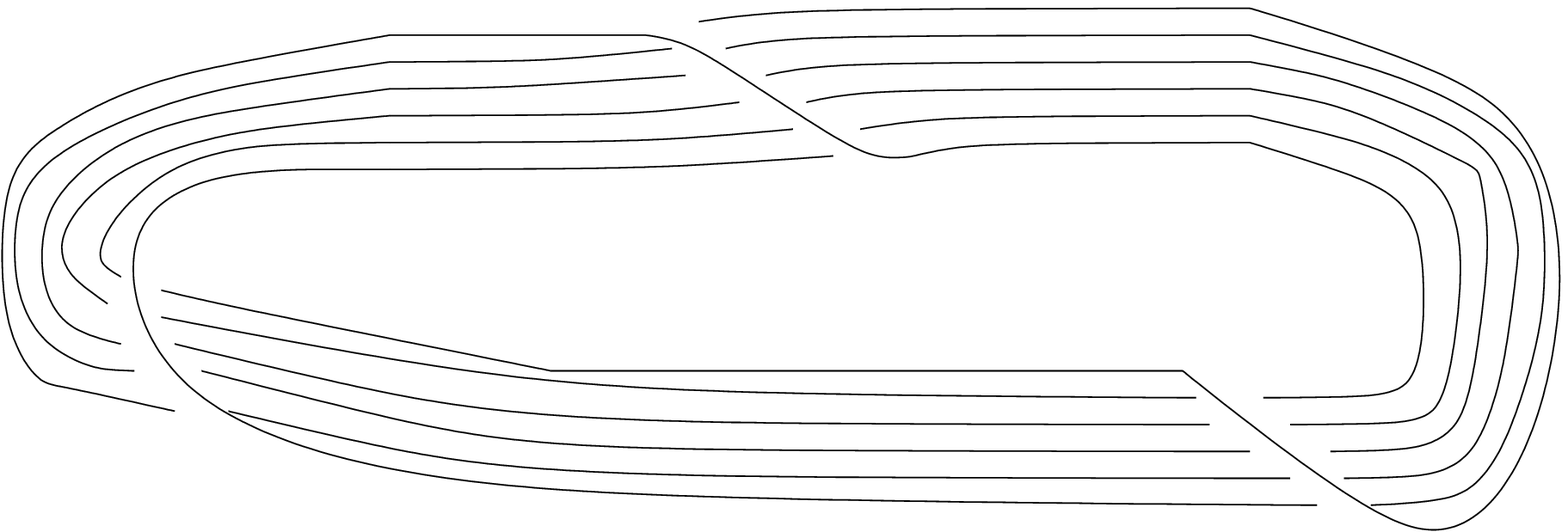,height=1.0in}}

\begin{center}
Fig. 1.1. \ The torus link of type $(3,6)$, $T_{3,6}$
\end{center}

\begin{example}\label{1.1}
Let $S^3 = \{ z_1,z_2\in C\times C:|z_1|^2+|z_2|^2 = 1\}$. 
Let us consider an action of $Z_p\oplus Z_q$ on $S^3$ 
which is generated by  $T_p$ and $T_q$,
where $T_p(z_1,z_2) = (e^{2\pi i/p} z_1,z_2)$ and
$T_q(z_1,z_2) = (z_1, e^{2\pi i/q}z_2)$. 
Show that this action preserves torus link of type
$(p,q)$. 
This link can be described as the following set
$\{ (z_1,z_2)\in S^3:\ z_1 = e^{2\pi i(\frac{t}{p} +\frac{k}{p})}, z_2 =
e^{2\pi it/q}\}$, where $t$ is an arbitrary real number and $k$ 
is an arbitrary integer. \\
If $p$ is co-prime to $q$ then $T_{p,q}$ is a knot and can be
parameterized by: 
$$R\ni t\mapsto (e^{2\pi it/p},\ e^{2\pi it/q})\in S^3\subset C^2.$$  
In this case $Z_p\oplus Z_q = Z_{pq}$ with a generator $T=T_pT_q$.
\end{example}

Subsequently, we will focus on the action of a cyclic group 
$Z_n$. We will mainly consider the case of an action on $S^3$ 
with a circle of fixed points. 
The new link invariants, (see Section 2), provide efficient criteria 
for examining such actions. 
\begin{definition}\label{1.2} 
A link is called $n$-periodic if there exists an action of $Z_n$ 
on $S^3$ which preserves the link\footnote{More precisely we require an 
existence of an embedding of circles, realizing the link, which is 
equivariant under the group action.}
 and the set of fixed points of the action 
is a circle disjoint from the link. 
If, moreover, the link is oriented then we assume that 
a generator of $Z_n$ preserves the orientation of the link 
or changes it globally (that is on every component). 
\end{definition} 

\section{Polynomial invariants of links}
We will describe in this section polynomial invariants of links,
crucial in periodicity criteria. We start from the skein (Jones-Conway or 
Homflypt) polynomial, $P_L(v,z)\in Z[v^{\pm 1},z^{\pm 1}]$ (i.e. $P_L(v,z)$
is a Laurent polynomial in variables $v$ and $z$).

\begin{definition}
The skein polynomial invariant of oriented links in $S^3$ 
can be characterized by 
the recursive relation (skein relation): 
\begin{enumerate} 
\item 
[(i)] $v^{-1}P_{L_+}(v,z)-vP_{L_-}(v,z)=zP_{L_0}(v,z)$,\ where $L_+,L_-$ and 
$L_0$ are three oriented link diagrams, which are the same outside a small 
disk in which they look as in Fig. 2.1, 
\end{enumerate} 
and the initial condition 
\begin{enumerate} 
\item 
[(ii)] $P_{T_1}=1$, \ where $T_1$ denotes the trivial 
knot. 
\end{enumerate}
\end{definition}

\ \\ 
\centerline{\psfig{figure=L+L-L0.eps}} 
\begin{center} 
Fig. 2.1
\end{center} 
We need some elementary properties of the skein polynomial:
\begin{proposition}
\begin{enumerate}
\item[(a)](\cite{L-M}) $z^{com(L)-1}P_L(v,z) \in Z[v^{\pm 1},z]$ where 
$com(L)$ is the number of link components of $L$. That is we do not use
negative powers of $z$.
Furthermore the constant term, with respect to variable $z$, is non-zero.
\item[(b)] $P_L(v,z) - P_{T_{com(L)}}$ is divisible by
$(v^{-1}-v)^2 - z^2$. Here $T_i$ denotes the trivial link of $i$
components.
\item[(c)] Let $\hat P_L(v,z)= z^{com(L)-1}P_L(v,z)$, then
$\hat P_L(v,z)- (v^{-1}-v)^{com(L)-1}$ is divisible by
$(v^{-1}-v)^2 - z^2$
\item[(d)] $P_L(v,z) \equiv (v^3z)^j P_{T_k}(v,z) \mod (\frac{v^6-1}{v^2-1},
z^2+1)$\\ for some $j$ and $k$ where $j\equiv (k-com(L)) \mod 2$. In particular
$P_L \equiv \delta P_{\bar L} \mod (\frac{v^6-1}{v^2-1},z^2+1)$ where
$\delta$ equals to $1$ or $-1$ and $\bar L$ denotes the mirror image of $L$.
\end{enumerate}
\end{proposition}
\begin{proof}
\begin{enumerate} 
\item[(a)-(c)]
Initial conditions and recurrence relation give immediately that 
$\hat P_L(v,z) \in Z[v^{\pm 1},z]$. (c) holds for trivial links 
from the definition (the difference is zero). To make inductive step, 
first notice that $\hat P_L(v,z)$ satisfies the skein relation 
$v^{-1}\hat P_{L_+}(v,z)-v\hat P_{L_-}(v,z)=z^{2\epsilon}\hat P_{L_0}(v,z)$ 
where $\epsilon$ is equal to $0$ in the case of the selfrossing of $L_+$
and $\epsilon =1$ in the mixed crossing case. 
We can rewrite the skein relation so that the inductive step is 
almost obvious:
$$v^{-1}(\hat P_{L_+} - (v^{-1} - v)^{com{L_+}-1}) 
-v(\hat P_{L_-} - (v^{-1} - v)^{com{L_-}-1}) =$$
$$ z^{2\epsilon}
(\hat P_{L_0} - (v^{-1} - v)^{com{L_0}-1}) + (v^{-1} - v)^{com{L_+}-2\epsilon}
(v^{-1} - v)^{2\epsilon}-z^{2\epsilon})$$
Namely the last term in the equality is always divisible 
by $(v^{-1}-v)^2 - z^2$ therefore if two other are divisible then the 
last one is.
We completed the proof of (c). From (c) follows that $\hat P_{L}(1,0)=1$,
thus the second part of (a) follows. (b) is just a weaker version of (c).
\item [(d)] It can be derived from the relation of $P_L(e^{2s\pi i/6},\pm 1)$
and the first homology group (modulo 3) of the 2-fold branched cyclic cover
of $(S^3,L)$ \cite{L-M,Yo-1,Ve}.
\end{enumerate}
\end{proof} 
\begin{definition}
Let $L$ be an unoriented diagram of a link. Then the Kauffman bracket 
polynomial $\langle L\rangle \in Z[A^{\mp 1}]$ is defined by the following
properties:
\begin{enumerate}
\item $\langle \bigcirc\rangle = 1$

\item $\langle\bigcirc\sqcup L\rangle = -(A^2+A^{-2})\langle L\rangle$

\item $\langle$ \parbox{1.2cm}{\psfig{figure=L+nmaly.eps}}$\rangle =
A\langle \parbox{0.9cm}{\psfig{figure=L0nmaly.eps}}
\rangle  
+
A^{-1}\langle
\parbox{0.9cm}{\psfig{figure=Linftynmaly.eps}}
\rangle$
\end{enumerate}
\end{definition}

The Kauffman bracket polynomial is
a variant of the Jones polynomial for oriented links. Namely,
for $A = t^{-\frac{1}{4}}$ and $D$ being an oriented diagram of $L$ we have
\begin{equation}
V_L(t) = (-A^3)^{-w(D)} <D>
\end{equation}
where $w(D)$ is the {\em planar writhe} ({\em twist} or {\em Tait number})
of $D$ equal to the algebraic sum of signs of crossings.

Kauffman noted that $V_L(t) = P_L(t,{\sqrt t} - \frac{1}{{\sqrt t}})$.
In particular 
$t^{-1}V_{L_+} -t V_{L_-} = ({\sqrt t} - \frac{1}{{\sqrt t}})V_{L_0}$.
Proposition 2.2 gives us:
\begin{corollary}
\begin{enumerate}
\item[(a)]\cite {Jo}.\ \
For a knot $K$, $V_K(t) \in Z[t^{\pm 1}]$ and $V_K(t) -1$ is divisible
by $(t-1)(t^3-1)$.
\item[(b)]\cite{Yo-1}.\ \
For a knot $K$, $V_K(t) -\delta V_K(t^{-1})$ is divisible by 
$\frac{t^3+1}{t+1}$, where $\delta$ equals to $1$ or $-1$.
\end{enumerate}
\end{corollary}

In the summer of 1985 (two weeks before discovering the ``bracket"),   L.~Kauffman invented
 another invariant of links \cite{Ka},\
 $F_L(a,  z) \in Z[a^{ \pm  1}, z^{ \pm  1}]$,
generalizing the polynomial discovered at the
beginning of 1985 by Brandt,   Lickorish,  Millett and Ho
\cite{B-L-M,Ho}. To define the
 Kauffman polynomial we first introduce the polynomial
invariant of link diagrams $\Lambda _D (a, z)$.
It is defined recursively by:
\begin{description}
\item
[(i)] $\Lambda _o (a,  z) = 1$,
\item [(ii)] $\Lambda_{{\psfig{figure=R+maly.eps}}}
 (a,  z) = a \Lambda_| (a, z);
\ \Lambda_{{\psfig{figure=R-maly.eps}}}
(a, z) = a^{-1} \Lambda_| (a,  z)$,
\item [(iii)]
$\Lambda_{D_+}(a,  z) + \Lambda_{D_-}(a, z) = z(\Lambda_{_0}(a, z) +
\Lambda_{D_\infty}(a, z))$.
\end{description}
 
The Kauffman polynomial of oriented links is defined by
$$F_L(a,  z) = a^{-w(D)} \Lambda _D(a, z)$$
where $D$ is any diagram of an oriented link $L$.
 
\begin{remark}\label{2.5}
Let $F^*(a,z)$ be the Dubrovnik variant of the Kauffman polynomial 
\footnote{Kauffman described the polynomial $F^*$ on a
postcard to Lickorish sent from Dubrovnik in September
'85.}. 
The polynomial $F^*$ satisfies initial conditions $F^*_{T_n}=(\frac{a-a^{-1}}
{z}+1)^{n-1}$ and the recursive relation $a^{w( D_+)}F^*_{D_+}-
a^{w(D_-)}F^*_{D_-} = z(a^{w( D_0)}F^*_{D_0}-a^{w( D_{\infty})}
F^*_{D_{\infty}})$.
Lickorish noted that the Dubrovnik polynomial is just a variant of the  
Kauffman polynomial:\\ 
$F^*_L(a,z)= (-1)^{(com(L)-1)}F_L(ia,-iz)$, (\cite{Li}).\\
\ \\
\ \\
Proof.\\
 We check the formula for trivial links and then show that if it
holds for three terms of the skein relation it also holds for the fourth one.
\begin{enumerate}
\item[(i)] $F_{T_n}(ia,-iz)=(-1)^{(n-1)} (\frac{a-a^{-1}}{-z}-1)^{n-1}=F^*_{T_n}(a,z)$.
\item[(ii)] $(ia)^{w( D_+)}F_{D_+}(ia,-iz)+ (ia)^{w(D_-)}F_{D_-}(ia,-iz) 
= (-iz)((ia)^{w( D_0)}F_{D_0}+(ia)^{w( D_{\infty})} 
F_{D_{\infty}})$. This reduces to $a^{w( D_+)}F^*_{D_+}-a^{w(D_-)}F^*_{D_-} = 
z(a^{w( D_0)}F^*_{D_0}+ (-1)^{com(D_0) - com(D_{\infty})}
a^{w( D_{\infty})}F^*{D_{\infty}}$,
which gives the skein relation for $F^*_{D}(a,z)$.
\end{enumerate} 
\end{remark}

\section{Periodic links and the Jones polynomial.}

Periodicity of links is reflected in the structure of
new polynomials of links.
We will describe this with details in the case of the Jones polynomial
(and the Kauffman bracket), mostly following \cite{Mu-2,T-1,P-3,Yo-1}
but also proving new results. Let $D$ be an r-periodic diagram of an 
unoriented r-periodic link, that is $\varphi (D)= D$ where $\varphi$ denote 
the rotation of $R^3$ along the vertical axis by the angle $2\pi/r$. 
In the coordinates of $R^3$ given by a complex number (for the first 
two real co-ordinates) and a real number 
(for the third) one gets: $\phi (z,t)= (e^{2\pi i/r}z,t)$. $\varphi$ is a
generator of the group $Z_r$ acting on $R^3$ (and $S^3=R^3 \cup \infty$).
In particular ${\varphi}^r=Id$. See Fig. 3.1.\\
\ \\ 
\centerline{\psfig{figure=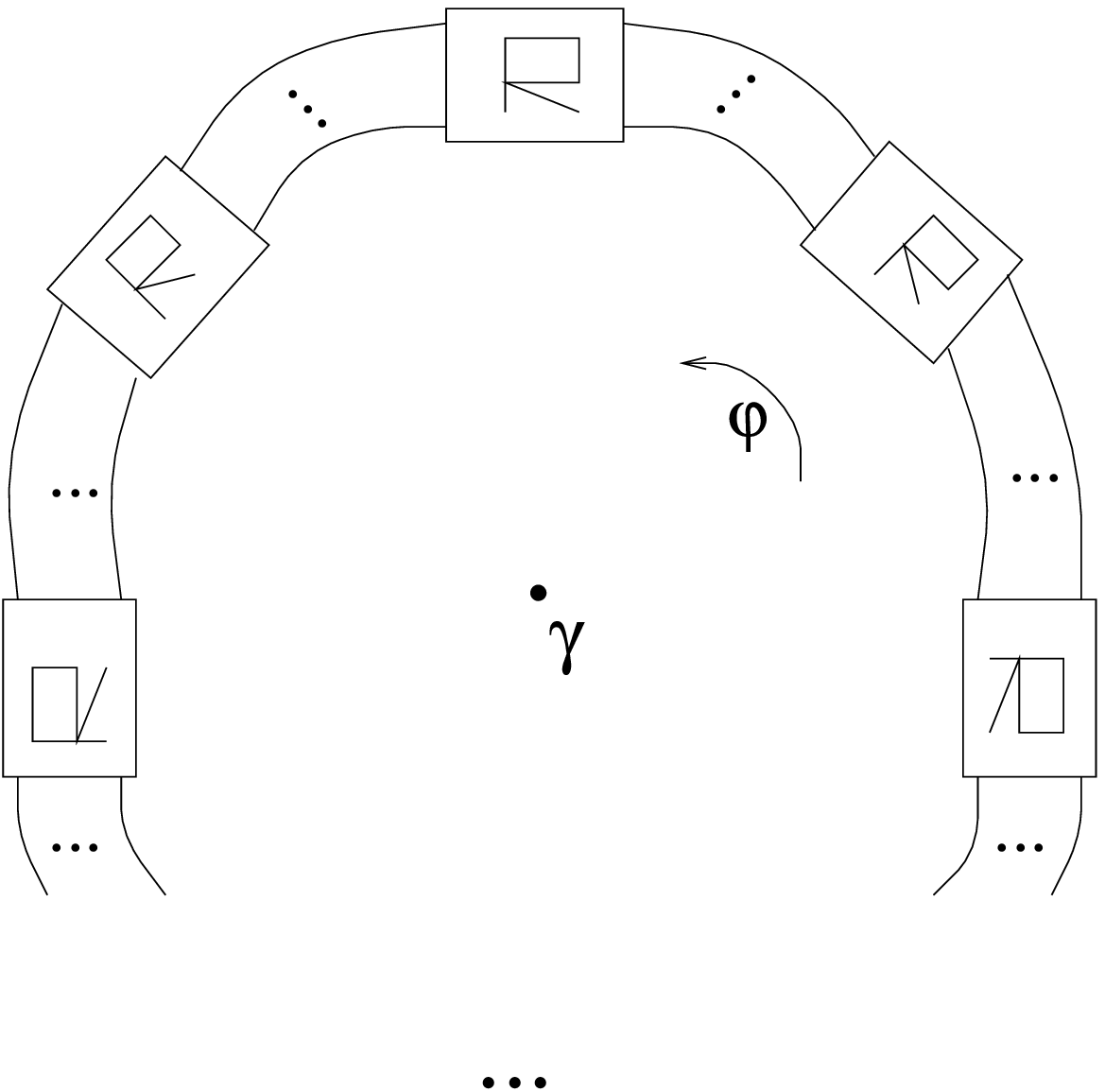,height=7.3cm}} 
\begin{center} 
Fig. 3.1 
\end{center} 
\ \\ 
By the positive solution to the Smith Conjecture \cite{Sm,Thur},
every r-periodic link has an r-periodic diagram, thus we can restrict
ourselves to considerations of these, easy to grasp, diagrams. With
the help of elementary group theory we have the following
fundamental lemma.
\begin{lemma} 
Let $D$ be an unoriented $r$-periodic link diagram, $r$ prime. 
Then the Kauffman bracket polynomial satisfies the following ``periodic"
formula:
$$ D_{sym ({\psfig{figure=L+nmaly.eps}})} \equiv 
A^r D_{sym ({\psfig{figure=L0nmaly.eps}})} + 
A^{-r} D_{sym ({\psfig{figure=Linftynmaly.eps}})} \mod r$$ 
where 
$D_{sym( 
{\psfig{figure=L+nmaly.eps}}) 
}, D_{sym({\psfig{figure=L0nmaly.eps}} 
)}$ and $D_{sym({\psfig{figure=Linftynmaly.eps}} 
)}$ 
denote three $\varphi$-invariant diagrams of links 
which are the same outside of the $Z_r$-orbit of a neighborhood of a 
fixed single crossing, $c$ (i.e. $c$, $\varphi (c)$,..., ${\varphi}^{r-1}(c)$)
 at which they differ by replacing 
{\psfig{figure=L+nmaly.eps}} 
by 
{\psfig{figure=L0nmaly.eps}} 
or 
{\psfig{figure=Linftynmaly.eps}}, 
respectively. 
\end{lemma}
\begin{proof}
Let us build the binary computational resolving tree of $D$ using
the Kauffman bracket skein relation for every crossing of the orbit
(under $Z_r$ action), $c$, $\varphi (c)$,..., ${\varphi}^{r-1}(c)$.
The tree has therefore $2^r$ leaves and a diagram at each leaf
contributes some value (polynomial) to the Kauffman bracket polynomial
of $D=D_{sym ({\psfig{figure=L+nmaly.eps}})}$, see Fig. 3.2. 
The idea of the proof is that only extreme
leaves, $D_{sym({\psfig{figure=L0nmaly.eps}})}$
 and $D_{sym({\psfig{figure=Linftynmaly.eps}} )}$, contribute
to $<D>$ $\mod r$. 

\ \\
\ \\ 
\centerline{\psfig{figure=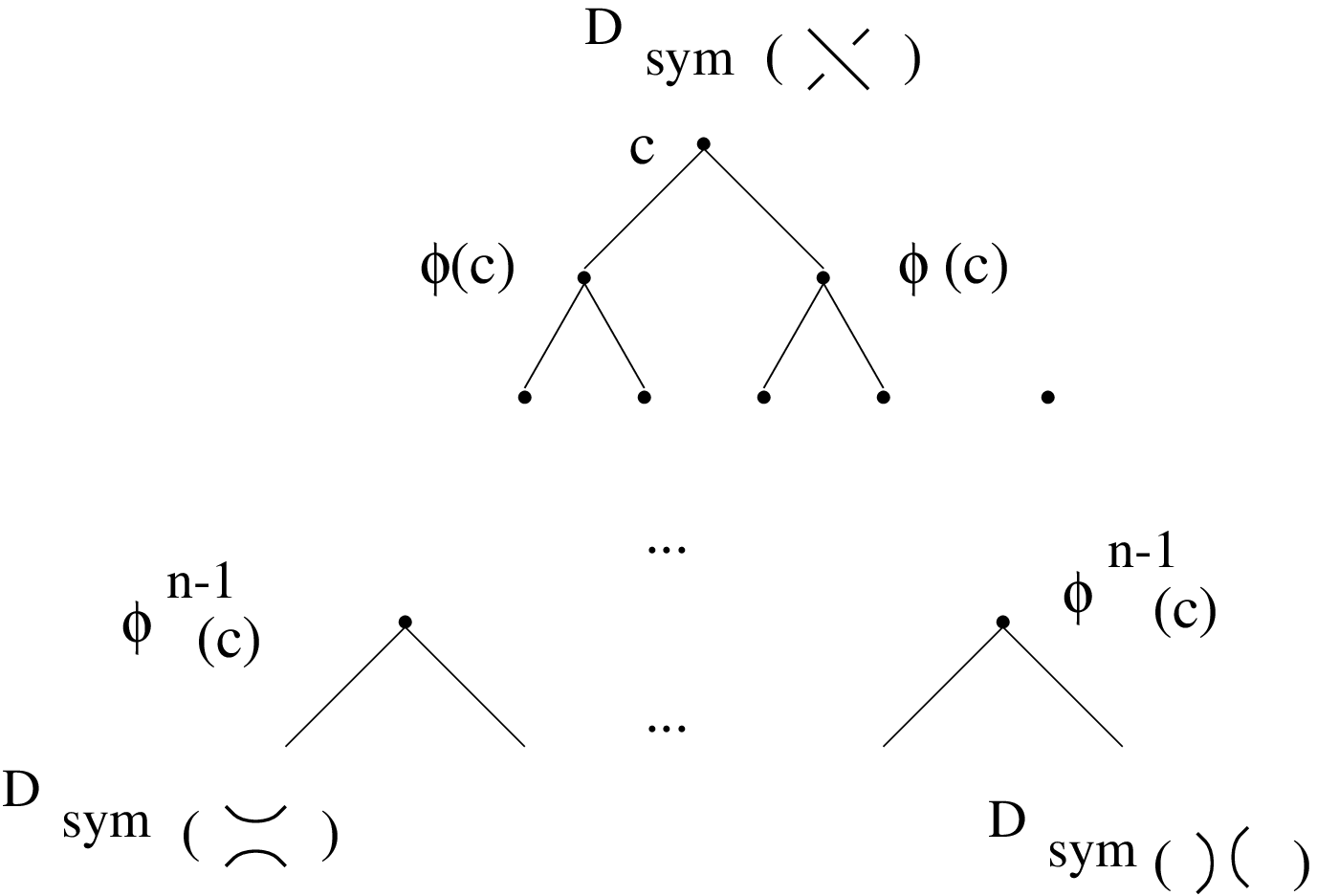,height=4.6cm}} 
\begin{center} 
Fig. 3.2 
\end{center} 
\ \\
We need now some, elementary, group theory:\\ 
Let a finite group $G$ acts on a set $X$, that is every element $g \in G$
``moves" $X$ ($g: X \to X$) and $hg(x) = h(g(x))$ for any $x\in X$ and
$h,g \in G$. Furthermore we require that the identity element of $G$ is not
moving $X$ ($e(x)=x$ where $e$ is the identity element of $G$). The orbit of
an element $x_0 $ in $X$ is the set of all elements of $X$ which can be
obtained from $x_0 $ by acting on it by $G$ (${\cal O}_{x_0}=
\{x\in X \ | \ \ there \ \ is \ \ g \ such  \ \ that \ \ g(x_0)=x \}$).
The standard, but important, fact of elementary group theory is that
the number of elements in an orbit divides the order (number of elements)
of the group. In particular if the group is equal to $Z_r$, $r$ prime,
then orbits of $Z_r$ action can have one element (such orbits are
called fixed points) or $r$ elements.\\
After this  long group theory digression, we go back to our leaf
diagrams of the binary computational resolving tree of $D$.
We claim that $Z_r$ acts on the leaf diagrams and the only fixed
points of the action are extreme leaves 
$D_{sym({\psfig{figure=L0nmaly.eps}})}$ and 
$D_{sym({\psfig{figure=Linftynmaly.eps}})}$. 
All other orbits have $r$ elements and they would cancel their contribution
to $<D>$ modulo $r$. To see these we introduce an adequate notation:
Let $c_i= \varphi (c)$ and $D^{c_0,...,c_{r-1}}_{s_0,...,s_{r-1}}$,
where $s_i = {\psfig{figure=L0nmaly.eps}}$ or 
${\psfig{figure=Linftynmaly.eps}}$, denote the diagram of a link obtained
from $D$ by smoothing the orbit of crossings, $c_0, ... c_{r-1}$ according
to indices $s_i$. $D^{c_0,...,c_{r-1}}_{s_0,...,s_{r-1}}$ are leaves of
our binary computational resolving tree of $D$, and the $Z_r$ action
can be fully described by the action on the indices $s_0,...,s_{r-1}$.
Namely $\varphi (s_0,...s_{r-2},s_{r-1}) = (s_{r-1},s_0,...,s_{r-2})$.
From this description it is clear that the only fixed point sequences
are $({\psfig{figure=L0nmaly.eps}},...,{\psfig{figure=L0nmaly.eps}},
{\psfig{figure=L0nmaly.eps}})$ and $({\psfig{figure=Linftynmaly.eps}},...,
{\psfig{figure=Linftynmaly.eps}},{\psfig{figure=Linftynmaly.eps}})$ and
that diagrams from the given orbit represent equivalent (ambient
isotopic) links, they just differ by the rotation of $R^3$. From this follows
that the contribution to $<D>$ of an $r$ element orbit is equal to $0$
modulo $r$. Thus all
leaves, but the fixed points, do not contribute to $<D>$ modulo $r$. 
The contribution of the fixed point leaves is expressed in the formula
of Lemma 3.1.
  \end{proof} 

\begin{corollary} Let $D^o$ be an oriented $r$-periodic link diagram, $r$ prime,
and $D$ the same diagram, orientation forgotten.
 \begin{enumerate} 
\item [(i]) 
Then the Kauffman bracket polynomial satisfies the following  
formula:
$$A^r <D_{sym ({\psfig{figure=L+nmaly.eps}})}> -
  A^{-r}<D_{sym ({\psfig{figure=L-nmaly.eps}})}> \equiv$$
$$
(A^{2r}-A^{-2r})<D_{sym ({\psfig{figure=L0nmaly.eps}})}> \mbox{ mod }(r)$$
\item [(ii])
Let $f_{D^o}(A)= (-A^3)^{-w(D^o)}<D>$ then 
$$-A^{4r}f_{D^o_{sym (+)}} + A^{-4r}f_{D^o_{sym (-)}} \equiv
(A^{2r}-A^{-2r})f_{D^o_{sym (0)}} \mbox{ mod }(r)$$
Here $D^o_{sym (+)}=D^o_{sym( 
{\psfig{figure=L+maly.eps}})) },\ D^o_{sym (-)}= 
D_{sym({\psfig{figure=L-maly.eps}} )}$ and 
$D^o_{sym (0)}=D_{sym({\psfig{figure=L0maly.eps}})}$ 
denote three $\varphi$-invariant oriented diagrams of links 
which are the same outside of the $Z_r$-orbit of a 
fixed single crossing and which at a neighborhood of the crossing 
differ by replacing 
{\psfig{figure=L+maly.eps}} 
by 
{\psfig{figure=L-maly.eps}} 
or 
{\psfig{figure=L0maly.eps}}; 
\item [(iii]) $$t^{-r}V_{D^o_{sym (+)}} - t^{r}V_{D^o_{sym (-)}} \equiv 
(t^{r/2}-t^{-r/2})V_{D^o_{sym( 0)}} \mbox{ mod }(r)$$ 
\end{enumerate}
\end{corollary}
\begin{proof}
\begin{enumerate} 
\item [(i]) 
Use Lemma 3.1 for $D_{sym ({\psfig{figure=L+nmaly.eps}})}$ and
$D_{sym ({\psfig{figure=L-nmaly.eps}})}$ and reduce the term 
$D_{sym ({\psfig{figure=Linftynmaly.eps}})}$.
\item [(ii])
(i) can be written as:\\
$(-A^3)^{w(D^o_{sym (+)})}A^rf_{D^o_{sym (+)}} -
(-A^3)^{w(D^o_{sym (-)})}A^{-r}f_{D^o_{sym (-)}} \equiv 
(A^{2r}-A^{-2r}) (-A^3)^{w(D^o_{sym (0)})}f_{D^o_{sym (0)}} \mbox{ mod }(r)$ 
and using the equality $w(D^o_{sym (+)}) = 
w(D^o_{sym (-)}) + 2r = w(D^o_{sym (0)})+r$,
one gets the congruence (ii).
\item [(iii)] This follows from (ii) by putting $V_L(t) = f_{L}(A)$, for 
$t=A^{-4}$.

\end{enumerate}
\end{proof}

Lemma 3.1 and Corollary 3.2 have several nice applications:\\
 to symmetric knots, periodic 3-manifolds and to analysis  of connections 
between skein modules of the base and covering space in a covering
(see Section 7). 
Below are two elementary but illustrative applications to periodic links.

\begin{theorem}[\cite{T-1,P-3}] 
\begin{enumerate} 
\item[(i)] 
If $L$ is an $r$-periodic oriented link ($r$ is a prime), 
then its Jones polynomial satisfies the relation 
$$V_L(t)\equiv V_L(t^{-1})\mbox{ mod }(r,t^r-1).$$ 
\item[(ii)]  
Let us consider a polynomial 
$\hat{V}_L(t) = (t^{\frac{1}{2}})^{ - 3lk(L)} V_L(t)$ 
which is an invariant of an ambient isotopy 
of unoriented links. $lk(L)$ denotes here the global linking number of $L$,
any orientation of $L$ gives the same $\hat{V}_L(t)$.\\
If $L$ is an $r$-periodic  unoriented link ($r$ is an odd prime), then 
$\hat{V}_L(t)\equiv\hat{V}_L(t^{-1})\mbox{ mod}(r,t^r-1)$. 
\end{enumerate} 
\end{theorem} 

\begin{theorem}\label{Theorem 3.4} 
Let $K$ be an $r$-periodic knot with linking number $k$ with the fixed 
point set axis. Then 
$$t^{r(k-1)/2}V_K(t) \equiv \frac{(t^{(k+1)/2}-t^{(-k-1)/2})- 
(t^{(k-1)/2}- t^{(1-k)/2})}{t-t^{-1}}\mod (r,t^r-1).$$  In particular: 
\begin{description} 
\item 
[(a)] If $k$ is odd then 
 $V_K(t) \equiv \frac {t^{k/2}+t^{-k/2}}{t^{1/2}+t^{-1/2} }\equiv 
t^{(k-1)/2} - t^{(k-3)/2}+...-t^{(3-k)/2}+t^{(1-k)/2} \mod (r,t^{r}-1)$, 
\item 
[(b)] If $k$ is even then 
 $t^{r/2}V_K(t) \equiv 
 \frac 
 {t^{k/2}+t^{-k/2}-2(-1)^{k/2}} 
 {t^{1/2}+t^{-1/2}} 
 +(-1)^{k/2} 
 \frac 
 {t^r+1} 
 {t^{1/2}+t^{-1/2} } 
\mod (r,t^r-1)$. 
\end{description} 
\end{theorem} 
\begin{proof} 
\begin{enumerate} 

\item[3.3(i)] 
From Corollary 3.2(iii) follows that  
$V_{D^o_{sym (+)}} \equiv  V_{D^o_{sym (-)}} \mod (r,t^{r}-1)$.
On the other hand we can change $D^o$ to its mirror image $\bar D^o$ by 
a sequence of changes of type $D^o_{sym (+)} \leftrightarrow
D^o_{sym (-)}$. Therefore $V_{D^o}(t) \equiv  V_{\bar D^o}(t) 
\mod (r,t^{r}-1)$. Theorem 3.3(i) follows because $V_{\bar D^o}(t) =
V_{D^o}(t^{-1})$.
\item[3.3(ii)]   
The link $L$ may be oriented so that $\varphi$ (a generator of the $Z_r$ 
action) preserves the orientation of $L$
(first we orient $L_*=L/Z_r$ the quotient of $L$ under the group action
 and then we lift the orientation up to $L$). For such an oriented $L$ we 
get from (i): 
$$V_L(t)\equiv V_L(t^{-1})\mbox{ mod }(r,t^r-1)$$ 
and thus 
$$(t^{\frac{1}{2}})^{3lk L}\hat{V}_L(t)\equiv(t^{\frac{1}{2}})^{-3lk 
L} \hat{V}_L(t^{-1})\mbox{ mod }(r,t^r-1)$$ and consequently 
$$\hat{V}_L(t)\equiv t^{-3lk(L)}\hat{V}_L(t^{-1})\mbox{ mod 
}(r,t^r-1).$$ For $r>2$, $lk(L)\equiv 0 \mod r$ so Theorem 3.3(ii) follows.

\item[3.4] We proceed as in the proof of 3.3(i) except that instead of
aiming at mirror image we aim the appropriate torus knot. To see this
it is convenient to think of our link $L$ being in a solid torus and
simplify its quotient $L_*=L/Z_n$ in the solid torus (see Section 7).
The formula for the Jones polynomial of the torus knot, $T_{(r,k)}$, 
was found by Jones \cite{Jo}:
$$V_{T_{(r,k)}}(t) = t^{r(k-1)/2} 
 \frac 
 {t^{(k+1)/2}-t^{(-k-1)/2}-t^r(t^{(k-1)/2}-t^{(1-k)/2}} 
 {t-t^{-1}}.$$  In Section 7, we discuss elementary proof of the formula
and give a short proof of the generalization of the Jones formula to the
solid torus (modulo $r$) using Lemma 3.1.
\end{enumerate} 
\end{proof} 

Theorem 3.3 is strong enough to allow Traczyk 
\cite{T-1} to complete periodicity tables for knots 
up to 10 crossings (for $r>3$), that is, to decide whether 
the knot $10_{101}$ is 7-periodic. 

\begin{example}\label{3.5} 
The knot $10_{101}$ (Fig.3.3) is not $r$-periodic for $r\geq 5$. 

\ \\ 
\centerline{\psfig{figure=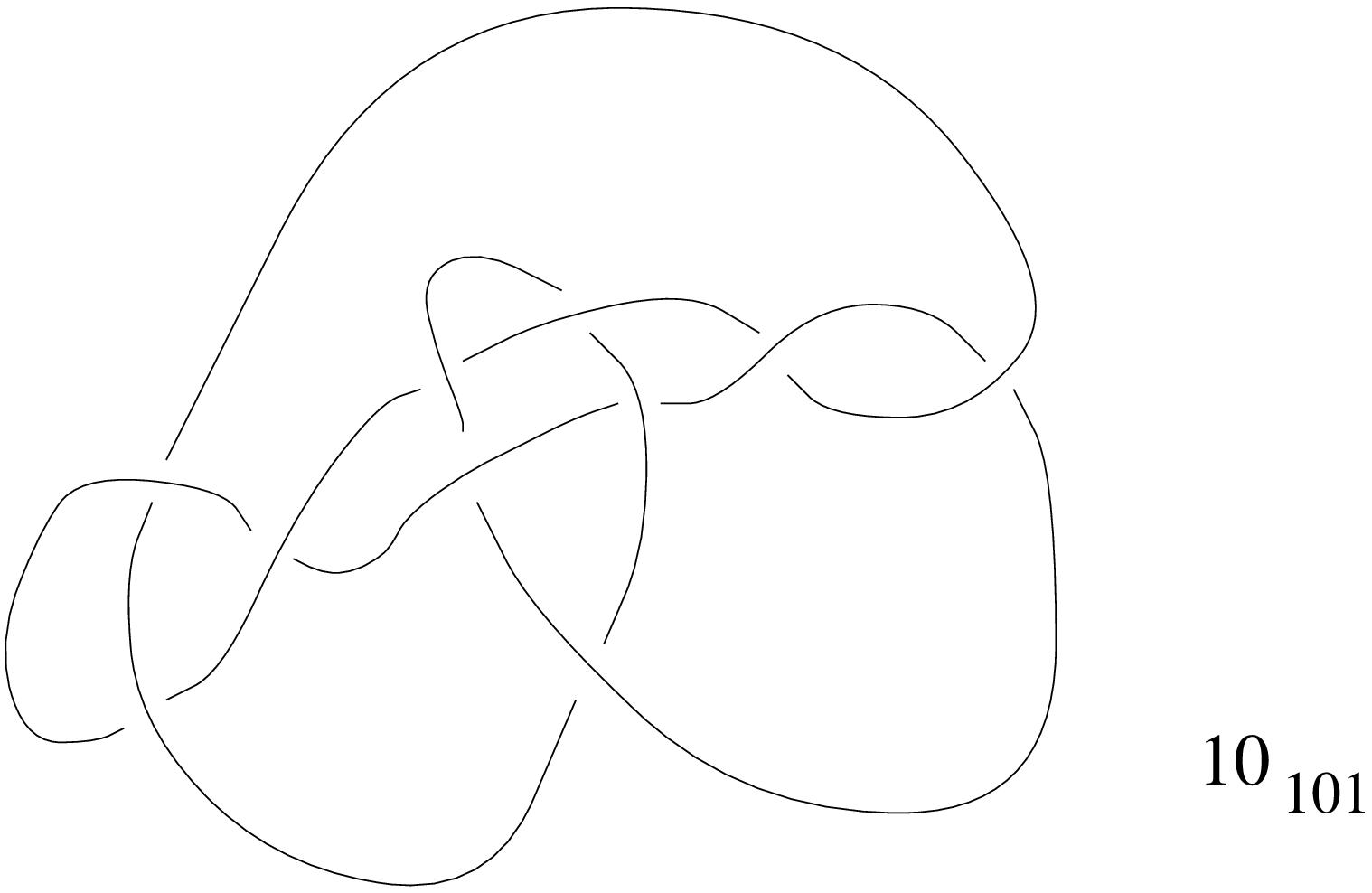,height=4.5cm}} 


\begin{center} 
Fig. 3.3
\end{center} 

The Jones polynomial for $10_{101}$ is equal to 
$$V_{10_{101}}(t) = t^2 - 3t^3 + 7t^4 - 
10t^5+14t^6-14t^7+13t^8-11t^9+7t^{10} - 4t^{11}+t^{12}.$$ 
Thus, for $r\geq 5$ it follows that $V_{10_{101}}(t)\not\equiv 
V_{10_{101}}(t^{-1}\mbox{ mod } (r,t^r-1)$. 
In particular, it follows that 
$V_{10_{101}}(t) \equiv t^{-1} + 2+3t^3\mbox{ mod }(5,t^5-1)$ and also 
$V_{10_{101}}(t) \equiv 3t^{-3}+5t^{-2}+6t+4t^2+4t^3\mbox{ mod 
}(7,t^7-1)$. 

\end{example} 
\begin{remark}
Theorem 3.3 and 3.4 do not work for knots and $r=3$. It is the case because
by Corollary 2.4, for any knot, $V_K(t) \equiv 1 \mod (t^3-1)$. One should
add that the classical Murasugi criterion using the Alexander polynomial,
Theorem 5.3 (\cite{Mu-1}) is working for $r=3$, and Traczyk developed
the method employing the skein polynomial, Theorem 4.10(a) (\cite{T-2,T-3}).
\end{remark}

Yokota proved in \cite{Yo-1} the following criterion for periodic knots 
which generalize Theorem 3.3 and is independent of Theorem 3.4.

\begin{theorem}\label{Theorem 3.7.} 
Let $K$ be an $r$-periodic knot ($r$ is an odd prime)
 with linking number $k$ with the fixed point set axis. Then 
\begin{enumerate} 
\item[(a)] 
 If $k$ is odd then 
    $V_K(t) \equiv V_K(t^{-1}) \mod (r,t^{2r}-1)$, 
\item [(b)] If $k$ is even then 
    $V_K(t) \equiv t^rV_K(t^{-1}) \mod (r,\frac{t^{2r}-1}{t+1})$. 
\end{enumerate} 
\end{theorem} 

We can extend Theorems 3.3, 3.4 and 3.7 (or rather show its limits) by
considering the following operations on link diagrams:
\begin{definition}\label{3.7}
\begin{enumerate}
\item [(i)] A $t_k$ move is an elementary operation on an oriented link
diagram $L$ resulting in the diagram $t_k(L)$ as shown on Fig. 3.4.
\ \\
\ \\
\centerline{\psfig{figure=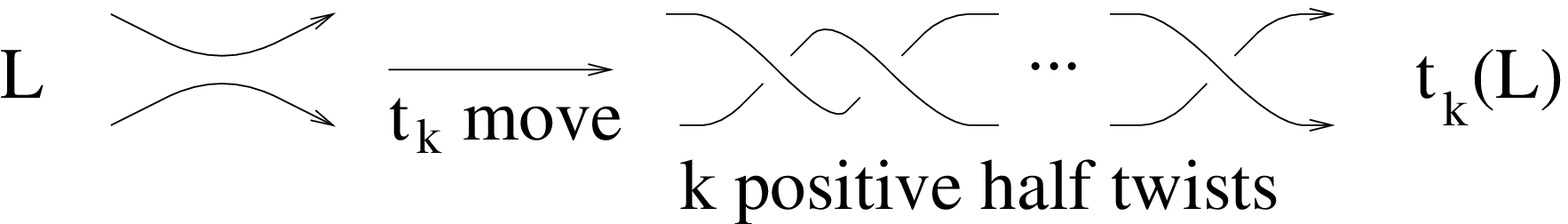,height=1.1cm}}  
\begin{center}  
Fig. 3.4 
\end{center} 
\item [(ii)] A $\bar t_k$ move, $k$ even, is an elementary operation on 
an oriented link diagram $L$ resulting in the diagram $t_k(L)$ 
as shown on Fig. 3.5.
\ \\
\ \\
\centerline{\psfig{figure=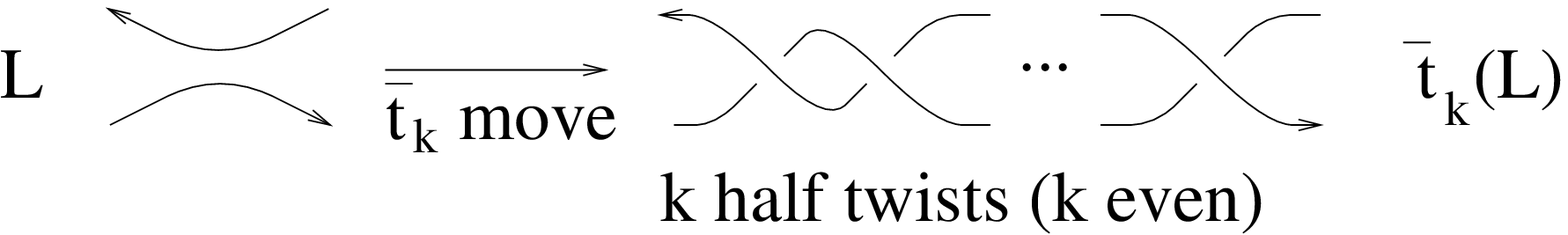,height=1.1cm}}  
\begin{center}   
Fig. 3.5  
\end{center}  
\item [(iii)] 
The local change in a link diagram which replaces parallel lines 
by $k$ positive half-twists is called a k-move; see Fig.3.6. 
\end{enumerate}
\end{definition}
\ \\
\centerline{\psfig{figure=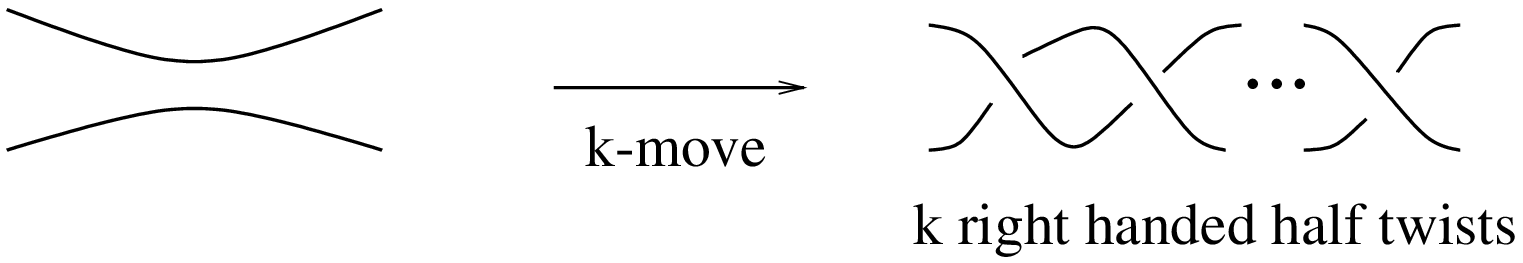,height=1.3cm}}
\begin{center} 
Fig. 3.6 
\end{center} 
\begin{lemma}
Let $L_k$ be an unoriented link obtained from $L_0$ by a $k$ move.
Then $<L_k> = A^k<L_0> + A^{-3k+2}\frac{A^{4k}-(-1)^k}{A^4 + 1}$.
\end{lemma}
\begin{proof}
It follows by an induction on $k$.
\end{proof}
\begin{corollary}\label{3.10}
If $L$ and $L'$ are oriented $t_{2r},\bar t_{2r}$ equivalent links 
(that is $L$ and $L'$ differ by a sequence of $t_{2r},\bar t_{2r}$-moves) 
then\\
 $V_L(t) \equiv t^{rj}V_{L'}(t) \mod (r, \frac{t^{2r}-1}{t+1})$, for
some integer $j$. 
\end{corollary}
\begin{proof}
For $k=2r$ we have from Lemma 3.9\  $<L_{2r}> = A^{2r}<L_0> + 
A^{-6r+2}(A^{4r}-1)\frac{A^{4r}+1}{A^4 + 1}<L_{\infty}>$. Thus 
$<L_{2r}> \equiv A^{2r}<L_0> \mod (A^{4r}-1)\frac{A^{4r}+1}{A^4 + 1}$ 
and  Corollary 3.10 easily follows.
\end{proof}

\section{Periodic links and the generalized Jones polynomials.}

We show here how periodicity of links is reflected in  regularities of
skein and Kauffman polynomials. We explore the same ideas which were
fundamental in Section 3, especially we use variations of Lemma 3.1.

Let ${\cal R}$ be a subring of the ring
$Z[v^{\mp 1},z^{\mp 1}]$ generated by $v^{\mp  1}$, $z$ and
$\frac{v^{-1}-v}{z}$. 
Let us note that  $z$ is not invertible in ${\cal R}$.

\begin{lemma}\label{4.1}
For any link $L$ its skein polynomial $P_L(v,z)$
is in the ring ${\cal R}$. 
\end{lemma}

\begin{proof}
For a trivial link $T_n$ with $n$ components we have
$P_{T_n}(v,z)=(\frac{v^{-1}-v}{z})^{n-1}\in {\cal R}$. 
Further, if
$P_{L_+}(v,z)$ (respectively $P_{L_-}(v,z)$) and $P_{L_0}(v,z)$ 
are in ${\cal R}$ then $P_{L_-}(v,z)$ (respectively $P_{L_+}(v,z)$) 
is in ${\cal R}$ as well.
This observation enables a standard induction to conclude \ref{4.1}.
Now we can formulate our criterion for $r$-periodic links. It has an 
especially simple form for a prime period (see Section 5 for a more general 
statement). 
\end{proof}

\begin{theorem}\label{4.2}
Let $L$ be an $r$-periodic oriented link and assume that $r$ is
a prime number. Then the skein polynomial $P_L(v,z)$
satisfies the relation
$$P_L(v,z)\equiv P_L(v^{-1},-z)\mbox{ mod }(r,z^r)$$
where $(r,z^r)$ is an ideal in ${\cal R}$ generated by $r$ and
$z^r$.
\end{theorem}

In order to apply Theorem \ref{4.2} effectively, we need the following fact.
\begin{lemma}\label{4.3}
Suppose that $w(v,z)\in{\cal R}$ is written in the form 
$w(v,z) = \sum_i u_i(v)z^i$, where
$u_i(v)\in Z[v^{\mp 1}]$. 
Then $w(v,z)\in(r,z^r)$ if and only if
for any $i\leq r$ the coefficient $u_i(v)$ is in the ideal
$(r,(v^{-1}-v)^{r-i})$.
\end{lemma}

\begin{proof} 
$\Leftarrow$\ Suppose $u_i(v) \in (r,(v^{-1}-v)^{r-i})$ for $i\leq r$.
Now $u_i(v)z^i \equiv (v^{-1}-v)^{r-i})z^ip(v) \equiv (\frac{v^{-1}-v}{z})
^{r-i}z^rp(v) \mod r$, where $p(v)\in Z[v^{\pm 1}]$.
Thus $u_i(v)z^i \in (r,z^r)$ and finally $w(v,z)\in(r,z^r)$.\\
$\Rightarrow$\ Suppose that $w(v,z)\in(r,z^r)$, that is, $w(v,z)\equiv
z^r\overline{w}(v,z)\mbox{mod }r$ for some
$\overline{w}(v,z)\in{\cal R}$. The element
$\overline{w}(v,z)$ can be uniquely written  as a sum 
$\overline{w}(v,z) =
z\overline{u}(v,z)+\sum_{j\geq 0}
(\frac{v^{-1}-v}{z})^j\overline{u}_j(v)$, where $\overline{u}(v,z)\in
Z[v^{\mp 1},z]$ and $\overline{u}_j(v)\in Z(v^{\mp 1})$. 
Thus, for 
$i\leq  r\ \ (j=r-i) $ 
we have $u_i(v)\equiv (v^{-1}-v)^{r-i}\overline{u}_{r-i}(v)\mbox{
mod } r$ and finally $u_i(v)\in (r,(v^{-1}-v)^{r-i})$ for $i\leq r$.
\end{proof}

\begin{example}\label{4.4}
Let us consider the knot $11_{388}$, in Perko's notation 
\cite{Per}, see Fig.4.1. 
The skein polynomial $P_{11_{388}}(v,z)$ is equal to
$$(3-5v^{-2}+4v^{-4}-v^{-6}) +(4-10v^{-2}+5v^{-4})z^2 +
(1-6v^{-2}+v^{-4})z^4 -v^{-2}z^6.$$ 
Let us consider the polynomial $P_{11_{388}}(v,z) - P_{11_{388}}(v^{-1},-z)$.
The coefficient $u_0(v)$ for this polynomial is equal to
$5(-v^{-2}+v^{2})+4(v^{-4}-v^{4})-v^{-6}+v^{6}$ and thus for $r\geq 7$ we have
$u_0(v)\not\in (r,(v^{-1}-v)^r) = (r,v^{-r}-v^{r})$. 
Now, from Lemma 4.3 we have $P_{11_{388}}(v,z) - P_{11_{388}}(v^{-1},-z)\not\in
(r,z^r)$. Therefore from Theorem 4.2 it follows that
the knot $11_{388}$ is not $r$-periodic for $r\geq 7$.

\ \\
\centerline{\psfig{figure=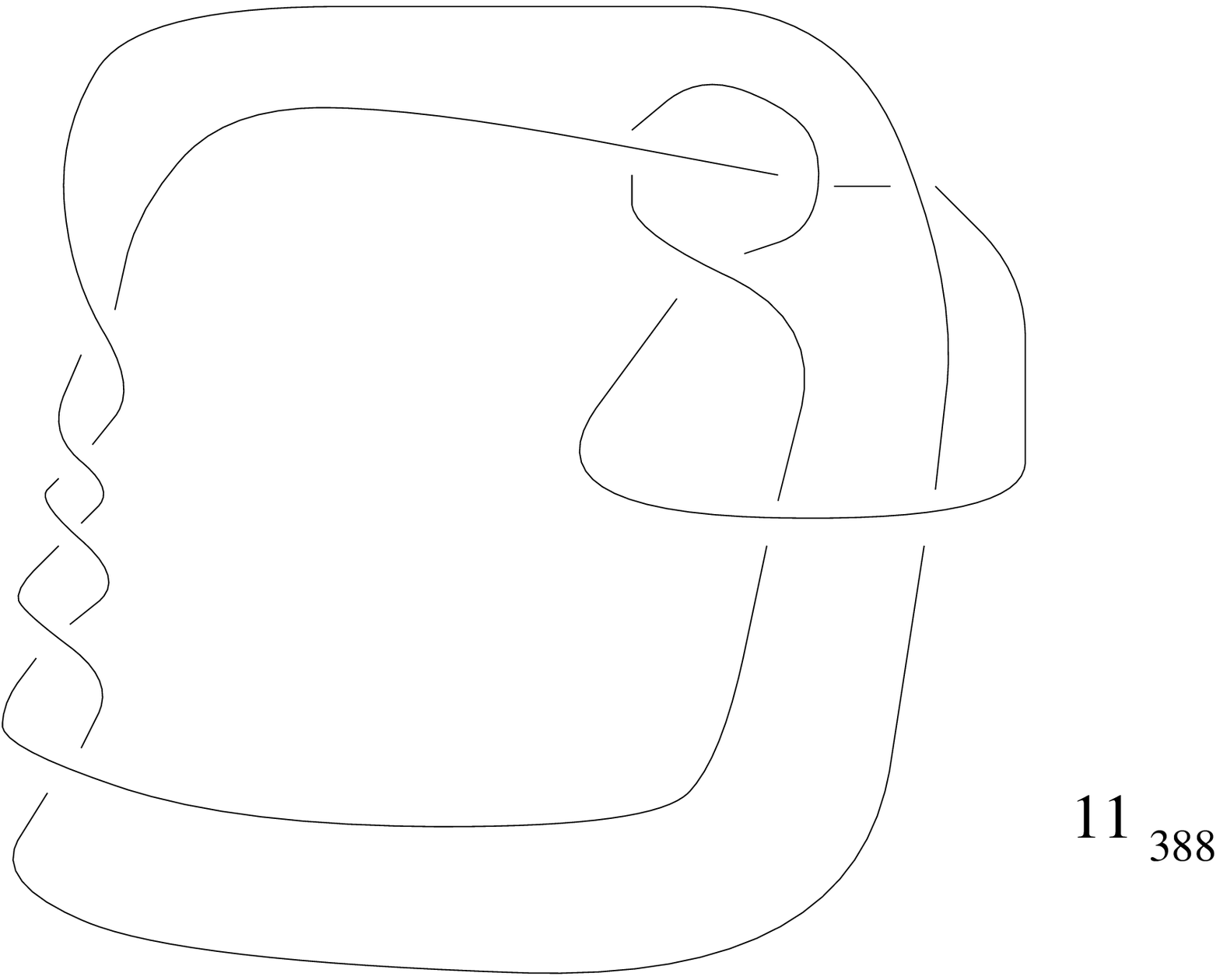,height=4.2cm}}
\begin{center}
Fig. 4.1
\end{center}

\end{example}

Theorem 3.3 is a corollary of Theorem 4.2 as 
$V_L(t)=P_L(t,t^{1/2}-t^{-1/2})$.

The periodicity criterion from Theorem 3.3 
is weaker than the one from Theorem 4.2: 
the knot $11_{388}$ from Example 4.4 has a symmetric
Jones polynomial
$$V_{11_{388}}(t) = V_{11_{388}}(t^{-1}) = t^{-2} - t^{-1} +
1-t+t^2,$$ 
and therefore Theorem 3.3 can not be applied in this case.
Theorem 3.4 is also not sufficient in this case (linking number $5$ cannot
be excluded as $V_{11_{388}}(t) = \frac{t^{5/2}+t^{-5/2}}{t^{1/2}+t^{-1/2}}
\equiv V_{T_{(r,5)}} \mod (t^r-1)$). 
\ \\
\ \\

Now let us consider the Kauffman polynomial $F_L(a,z)$. 
Let ${\cal R}'$ be a subring of of $Z[a^{\mp 1},z^{\mp 1}]$
generated by $a^{\mp 1},z$ and $\frac{a+a^{-1}}{z}$. It is easy to check
that Kauffman polynomials of links are in ${\cal R}'$ (compare Lemma 4.3).

\begin{theorem}\label{4.5}
Let $L$ be an $r$-periodic oriented link and let $r$ be  a prime number.
Then the Kauffman polynomial $F_L$ satisfies the following relation
$$F_L(a,z)\equiv F_L(a^{-1},z)\mbox{ mod }(r,z^r),$$
where $(r,z^r)$ is the ideal in ${\cal R}'$ generated by $r$ and
$z^r$.
\end{theorem}

In order to apply Theorem \ref{4.5} we will use
the appropriate version of Lemma 4.3.
\begin{lemma}
Suppose that $w(a,z)\in{\cal R}'$ is written in the form 
$w(a,z) = \sum_i v_i(a)z^i$, where 
$v_i(a)\in Z[a^{\mp 1}]$. 
Then $w(a,z)\in(r,z^r)$ if and only if 
for any $i\leq r$ the coefficient $v_i(a)$ is in the ideal 
$(r,(a+a^{-1})^{r-i})$. 
\end{lemma}

\begin{example}\label{4.7}
Let us consider the knot $10_{48}$ from Rolfsen's book \cite{Ro},
see Fig.4.2. 
This knot has a symmetric skein polynomial, that is 
$P_{10_{48}}(v,z) = P_{10_{48}}(v^{-1},-z)$. 
Consequently, Theorem 4.2
 can not be applied to examine the periodicity of this knot.
So let us apply the Kauffman polynomial to show that
the knot $10_{48}$ is not $r$-periodic for $r\geq 7$. 
It can be calculated (\cite{D-T,P-1})
that $F_{10_{48}}(a,z) - F_{10_{48}}(a^{-1},z) =
z(a^5+3a^3+2a-2a^{-1}-3a^{-3}-a^{-5})+z^2(\ldots)$.

Now let us apply Lemma 4.6 for $i=1$ and let us note that for
$r\geq 7$ we have $a^5+3a^3+2a-2a^{-1}-3a^{-3}-a^{-5}\not\in
(r,(a+a^{-1})^{r-1})$ 
(Note that for $r=5$ we have $a^5+3a^3+2a-2a^{-1}-3a^{-3}-a^{-5}
= a(a+a^{-1})^4-a^{-1}(a+a^{-1})^4 = (a-a^{-1})(a+a^{-1})^4\in
(5,(a+a^{-1})^4)$).

\end{example}
\ \\ 
\centerline{\psfig{figure=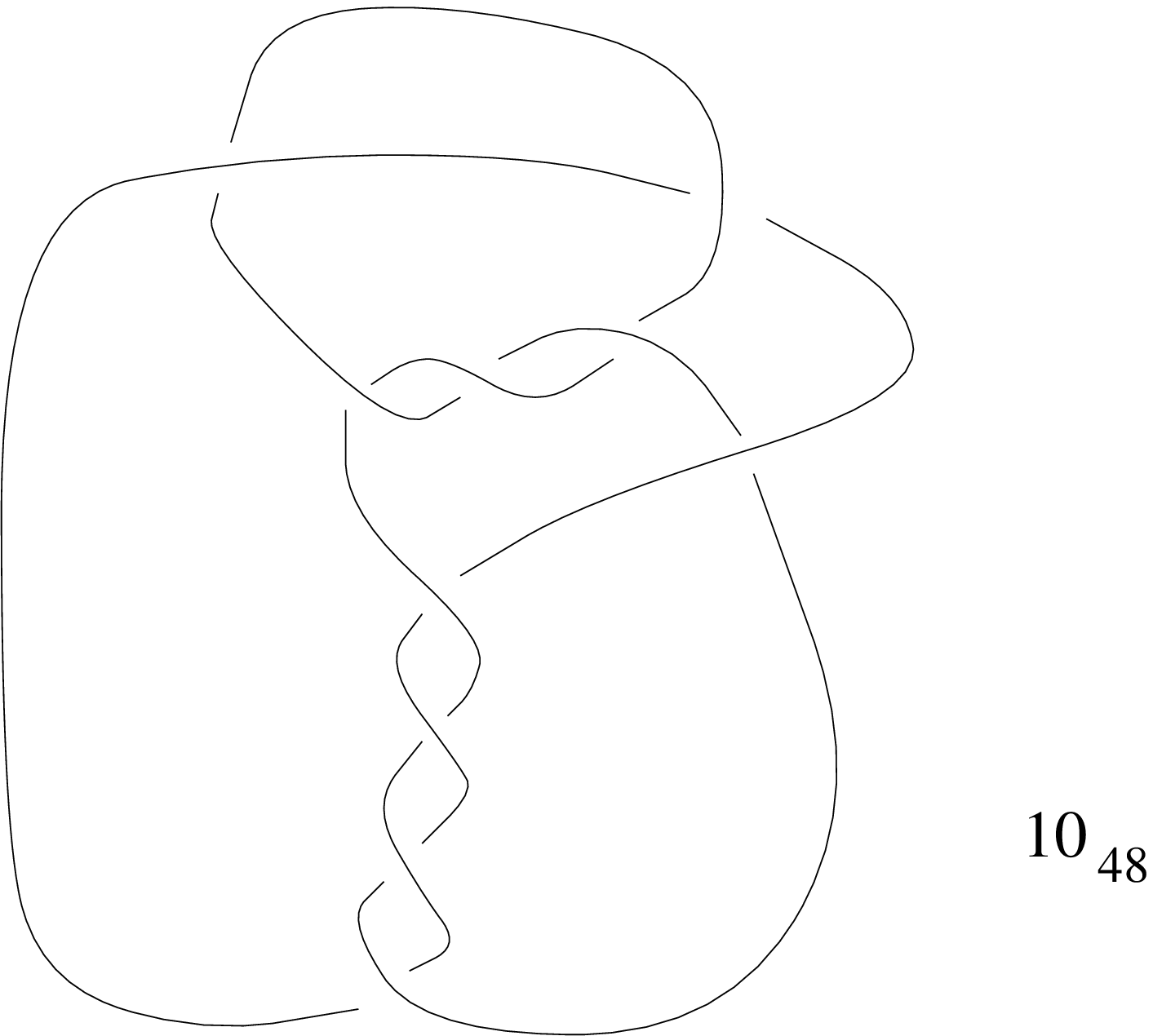,height=4.5cm}}
\begin{center}
Fig. 4.2
\end{center}

We can reformulate Theorem 4.5 in terms of the Dubrovnik version of the
Kauffman polynomial:

\begin{theorem}\label{4.8} 
Let $L$ be an $r$-periodic oriented link and let $r$ be  a prime number.
Then the Dubrovnik polynomial 
$F^*_L\in {\cal R}''=Z[a^{\pm 1},\frac{a-a^{-1}}{z},z]\subset 
Z[a^{\pm 1}, z^{\pm 1}] $ 
satisfies the following relation 
$$F^*_L(a,z)\equiv F_L^*(a^{-1},-z)\mbox{ mod }(r,z^r),$$ 
where $(r,z^r)$ is the ideal in ${\cal R}''$ generated by $r$ and 
$z^r$. 
\end{theorem} 

In the last part of this section
we will show how to strengthen Theorems 4.2 and 4.5.

One can also modify our method so it applies to symmetric links with one
fixed point component. We consider more general setting in Section 7
(symmetric links in the solid torus).

Proof of Theorems 4.3 and 4.5 is very similar to that of Theorem 3.3.
Instead of Lemma 3.1 we use the following main lemma (which is of
independent interest), proof of which is again following the same
principle ($Z_r$ acting on the leaves of a computational tree) as the
proof of Lemma 3.1.

\begin{lemma}\label{4.9}
\begin{enumerate}
\item[(i)]
Let     $L_{sym({\psfig{figure=L+maly.eps}})
},       L_{sym({\psfig{figure=L-maly.eps}} 
)}$ and $L_{sym({\psfig{figure=L0maly.eps}} 
)}$ 
denote three $\varphi$-invariant diagrams of links
which are the same outside of the $Z_r$-orbit of a 
fixed single crossing and which at the crossing
differ by replacing 
{\psfig{figure=L+maly.eps}}
by 
{\psfig{figure=L-maly.eps}}
or
{\psfig{figure=L0maly.eps}},
respectively. $\varphi$, as before denotes the generation of the $Z_r$ action.
Then 
$$a^r P_{L_{sym({\psfig{figure=L+maly.eps}}
)}}(a,z) + a^{-r}P_{L_{sym({\psfig{figure=L-maly.eps}}
)}}(a,z) = z^rP_{L_{sym({\psfig{figure=L0maly.eps}}
)}}(a,z) \mod r.$$
\item[(ii)]
Consider four $r$-periodic unoriented diagrams 
$L_{sym({\psfig{figure=L+nmaly.eps}})} , 
L_{sym({\psfig{figure=L-nmaly.eps}}) }, 
L_{ 
sym{\psfig{figure=L0nmaly.eps}}) 
} $ 
and $L_{sym({\psfig{figure=Linftynmaly.eps}})} $. Than the Kauffman
polynomial of unoriented diagrams diagrams, $\Lambda_L(a,z)$, 
satisfies:
$$\Lambda_{sym({\psfig{figure=L+nmaly.eps}}) }  + 
\Lambda_{sym({\psfig{figure=L-nmaly.eps}}) }  \equiv 
z^r(\Lambda_{sym({\psfig{figure=L0nmaly.eps}}) }  +
\Lambda_{sym({\psfig{figure=Linftynmaly.eps}}) } ) \mod r.$$

\end{enumerate}
\end{lemma}

Traczyk \cite{T-3} and Yokota \cite{Yo-2} substantially generalized Theorems
4.2 and 4.5. Lemma 4.9 is still crucial in their proofs, but detailed
study of the skein polynomial of the torus knots is also needed.
I simplified their proof by using Jaeger composition product \cite{P-5}.

\begin{theorem}\label{Theorem 4.10} 
Let K be an r-periodic knot ($r$ an odd prime number)
 with linking number with the rotation axis 
equal to k  and $P_K(v,z)=\sum_{i=0} P_{2i}z^{2i}$ then: 
\begin{description} 
\item 
[(a)] (Traczyk) If $P_0(K) = \sum a_{2i}v^{2i}$ then $a_{2i} \equiv a_{2i+2} 
\mod r$ except 
possibly when $2i+1 \equiv \pm k \mod r$. 
\item 
[(b)] (Yokota) $P_{2i}(K) \equiv b_{2i}P_0(K)\mod r$  for $2i < r-1$ , 
where numbers 
$ b_{2i} $ depends only on $r$ and $k \mod r$. 
\end{description} 
\end{theorem} 

One can also use the Jaeger's skein state model for Kauffman 
polynomial to give periodicity criteria yielded by the Kauffman polynomial.
I have never written details of the above idea (from June 1992),
 and Yokota proved independently criteria yielded by the Kauffman polynomial
\cite{Yo-3,Yo-4}.
\begin{theorem}[(Yokota)] Let K be an r-periodic knot ($r$ an odd prime number)
 with the linking number with the rotation axis equal to $k$, and let 
$(\frac{a-a^{-1}}{z}+1)F^*(a,z)=\sum_{i=0} F_{i}z^{i}$
($F_0(a,z)=P_0(v,z)$ for $v=a^{-1}$),
then for $2i \leq r-3$:
\begin{enumerate}
\item[(i)] $F _{2i}(K) \equiv b_{2i}F_0(K)\mod r$, \\
$F_{2i+1}\equiv 0 \mod r$ except $i=0$ where $F_1\in Z[a^{\pm r}] \mod r$ 
\item[(ii)] Let $P^*(a,z)= P(v,z)$ for $a=v^{-1}$ and
$J_K(a,z)=(\frac{a-a^{-1}}{z}+1)F^*_K(a,z) -(\frac{a-a^{-1}}{z})P^*_K(a,z)$
and let  $J_K(a,z)=\sum_{i=0} J_{i}z^{i}$.\ Then\\
for $0\leq i \leq r-1$ $J_i(a,k)\equiv 0 \mod r$ except for $i=1$ when
$J_1(a;K)\in Z[a^{\pm r}] \mod r$\\
Define $J_{i,l}(a;K)$ as a polynomial obtained by gathering all terms
in $J_i(a;K)$ which have degree $\pm l \mod r$. Then:\\
For each $l$ and for $0\leq 2i \leq r-3$
$$ J_{r+2i,l}(a;K) \equiv b_{i,k}J_{r,k}(a,K) \mod r$$
$ J_{r+2i+1}(a;K) \equiv 0 \mod r$ except $J_{r+1}(a;K)$ which modulo $r$
is in $Z_r[a^{\pm r}]$.
\end{enumerate}
\end{theorem}
\section{$r^q$-periodic links and Vassiliev invariants.}

The criteria of $r$-periodicity, which we have discussed 
before, can be partially extended to the case of
$r^q$-periodic links.
We assume that $r$ is a prime number and the fixed point set
of the action of $Z_{r^q}$ is a circle disjoint from the link in question
(trivial knot by the Smith Conjecture)

We will not repeat here all criteria where $r$ is generalized to $r^q$
\cite{P-3}, but instead we will list one pretty general criterion using
Vassiliev-Gusarov invariants (compare \cite{P-4}).
\begin{definition}
Let ${\cal K}^{sg}$ denote the set of singular oriented knots in $S^3$ where 
we allow only immersion of $S^1$ with, possibly, double points, up to 
ambient isotopy. Let $Z{\cal K}^{sg}$ denote the free abelian group 
generated by elements of ${\cal K}^{sg}$ (i.e. formal linear combinations of
singular links).
In the  group $Z{\cal K}^{sg}$ we 
consider resolving singularity relations $\sim$: $K_{cr} = K_+ - 
K_-$ ; see Fig. 5.1. 
\ \\ 
\centerline{\psfig{figure=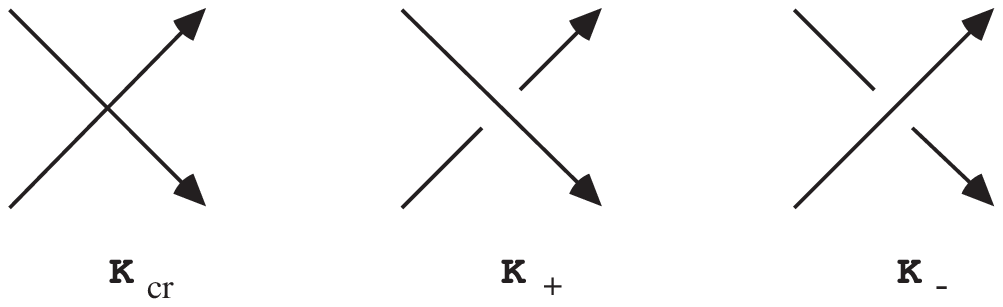}} 

\begin{center} 
Fig. 5.1 
\end{center} 

$Z{\cal K}^{sg}/\sim$ is clearly $Z$-isomorphic 
to $Z\cal K$. Let $C_m$ be a 
subgroup of $Z{\cal K}^{sg}/\sim =Z\cal K$ generated by immersed knots 
with $m$ double points. Let $A$ be any abelian group.
The $m$'th Vassiliev-Gusarov invariant is a homomorphism $f: Z{\cal K} \to A$
such that $f(C_{m+1}) =0$.
\end{definition}  
\begin{theorem}
\begin{enumerate} 
\item[(i)] Let $K$ be an oriented $r^q$-periodic knot and $f$ a 
Vassiliev-Gusarov invariant of degree $m<r^q$. Then
$f(K)\equiv f(\bar K) \mod r$, where $\bar K$ is the mirror image of $K$.
\item[(ii)] Let $K$ be an oriented $r^q$-periodic knot with the
 linking number $k$ with the fixed point set axis, and let $f$  be a 
Vassiliev-Gusarov invariant of degree $m<r^q$. Then 
$f(K)\equiv f(T_{(r^q,k})$ where $T_{(r^q,k)}$ is the torus knot of type
$(r^q,k)$.
\end{enumerate}
\end{theorem}

Proof of Theorem 5.2  is similar to the previous one and again bases
on the fundamental observation that 
$f(K_{sym(+)})-f(K_{sym(-)})\equiv 0 \mod r$.

Our method allows also to prove quickly the classical Murasugi congruence
for $r^q$ periodic knot, using the Alexander polynomial.
\begin{theorem}\label{5.3} 
Show that by applying our method to 
Alexander polynomial we obtain the following version of 
a theorem of Murasugi:
$${\Delta}_L(t)\equiv {\Delta}_{L_{*}}^{r^q}(t)
(1+t+t^2+...+t^{\lambda -1})^{r^q-1} \mod r $$ where $L_{*}$ 
is the quotient of an $r^q$-periodic link $L$ and 
$\lambda$ is the linking number of $L$ and $z$ axis.
\end{theorem}
\begin{proof} {\bf Sketch.}\\
Construct the binary computational  tree of the Alexander 
polynomial  of $L_*$ and the associated binary tree for Alexander polynomial
$L$ modulo $r$. For the Alexander polynomial of $L_*$ we use the skein 
relation $\Delta_{L_{*+}}(t) - \Delta_{L_{*-}}(t) = (t^{1/2} - 
t^{-1/2})\Delta_{L_{*0}}(t)$ and for the $r^q$ periodic link $L$ the congruence
(related to Lemma 4.8(i)):
$$\Delta_{L_{sym(+)}}(t) - \Delta_{L_{sym(-)}}(t) \equiv
(t^{r^q/2}-t^{-r^q/2})\Delta_{L_{sym(0)}}(t) \mod r .$$
Finally we use the fact that Alexander polynomial of a split
link is zero, and that $\Delta_{T_(r^q,\lambda)}\equiv 
(1+t+t^2+...+t^{\lambda -1})^{r^q-1} \mod r $. 
\end{proof}

\begin{remark}\label{5.9}
If $r^q=2$ then the formula from the previous exercise 
reduces to:
$${\Delta}_L(t)\equiv {\Delta}_{L_{*}}^2(t)
(1+t+t^2+...+t^{\lambda -1})\mod 2 .$$
Similar formula can be proven, for other $Z_2$-symmetry of links. Namely:\\
a knot (or an oriented link) in $R^3$ is 
called strongly plus amphicheiral if it has a realization in $R^3$ which 
is preserved by a (changing orientation) central symmetry 
$((x,y,z) \to (-x,-y,-z))$ ; ``plus" means that the involution is preserving 
orientation of the link. 
One can show,  using ``skein" considerations, as in the case 
of Theorem 5.3,  that if $L$ is a strongly $+$ amphicheiral link, then 
modulo 2 the polynomial ${\Delta}_L(t)$ is a square of another polynomial. 

Hartley and Kawauchi \cite{H-K} proved that ${\Delta}_L(t)$ is a square in
general. For example if $L(2m+1)$ is a Turks head link - the closure of 
the 3-string braid $(\sigma_1\sigma_2^{-1})^{2m+1}$, then its Alexander 
polynomial satisfies:
$${\Delta}_{L(2m+1)}= (\frac{a^{m+1}-a^{-m-1}}{a-a^{-1}})^2 $$
where $a+a^{-1}=1-t-t^{-1}$ and the Alexander polynomial is described
up to  an invertible element in $Z[t^{\pm 1}]$. This formula follows
immediately by considering the Burau representation of the 3-string
braid group \cite{Bi,Bu}.
\end{remark}

\section{Lissajous knots and billiard knots.}
A Lissajous knot $K$ is a knot in $R^3$ given by the parametric 
equations $$x = cos(\eta_x t + \phi_x )$$ $$y = cos(\eta_y t + \phi_y)$$ 
$$z = cos(\eta_z t + \phi_z)$$ for integers $\eta_x , \eta_y , \eta_z$. 
A {\em Lissajous link} is a collection of disjoint Lissajous knots. 

The fundamental question was asked in \cite{BHJS}: 
which knots are Lissajous?\\ 
It was shown in \cite{BHJS} and \cite{J-P} that a Lissajous knot is 
a $Z_2$-symmetric knot ($2$-periodic with a linking number with the axis equal
to $\pm 1$ or strongly plus amphicheiral) so a "random" knot is not
Lissajous (for example a nontrivial torus knot is not Lissajous).
Lamm constructed infinite family of different Lissajous knots \cite{La}. 

One defines a {\em billiard knot} (or racquetball knot) as the trajectory 
inside a cube of 
a ball which leaves a wall at rational angles with respect to the 
natural frame, and travels in a straight line except for reflecting 
perfectly off the walls; generically it will miss the corners and 
edges, and will form a knot. We show in \cite{J-P} that these knots 
are precisely the same as the Lissajous knots. 
We define general billiard knots, e.g. taking another polyhedron instead 
of the ball, considering a non-Euclidean metric, or considering the 
trajectory of a ball in the configuration space of a flat billiard. 
We will illustrate these by various examples. For instance, the trefoil knot 
is not a Lissajous knot but we can easily realize it as a billiard knot 
in a room with a regular triangular floor.

\ \\ 
\centerline{\psfig{figure=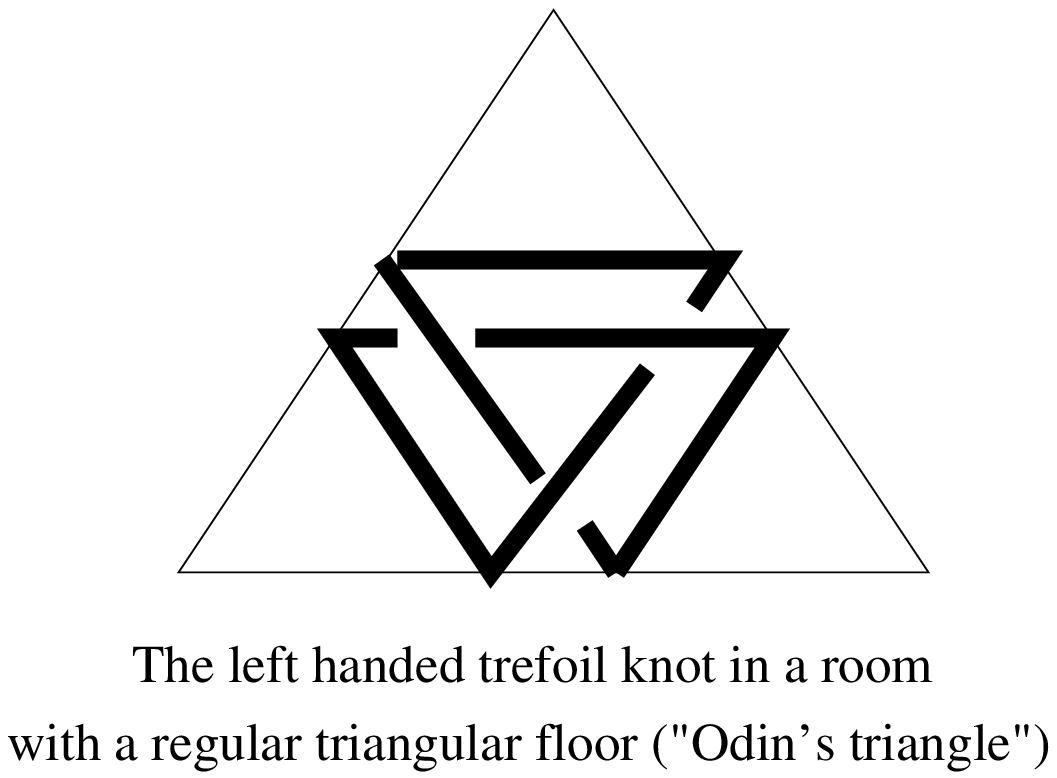,height=5.3cm}} 
\begin{center} 
Fig. 6.1 
\end{center} 
\begin{theorem}\label{6.1} 
Lissajous knots and billiard knots in a cube are the same up to ambient 
isotopy. 
\end{theorem}
A billiard knot (or link), is a simple closed trajectory (trajectories) of a 
ball in a 3-dimensional billiard table. 
The simplest billiards to consider would be polytope (finite convex 
polyhedra in $R^3$). But even for Platonian bodies we know nothing of 
the knots they support except in the case of the cube. It seems that 
polytopes which are the products of  polygons 
and the interval ($[-1,1]$) (i.e. polygonal prisms) are more accessible. 
This is the case because diagrams of knots are billiard trajectories in 
2-dimensional tables. We will list some examples below (compare \cite{Ta}). 

\begin{example}\label{6.2}
\begin{enumerate} 
\item [(i)] 
The trivial knot and the trefoil knot are the trajectories of a ball in 
a room (prism) with an acute triangular floor. 
In Fig.6.2(a), the diagram of the 
trivial knot is an inscribed triangle $\Delta_I$ whose vertices are the feet of 
the triangle's altitudes. If we move the first vertex of $\Delta_I$ slightly, 
each of its edges splits into two and we get the diagram of the trefoil. We 
should be careful with the altitude of the trajectory: We start from
level 1 at the vertex close to the vertex of $\Delta_I$ and opposite to 
the shortest edge of $\Delta_I$. Then we choose the vertical parameter so 
that the trajectory has 3 maxima and three minima (Fig.6.2(b)). 

\ \\ 
\centerline{\psfig{figure=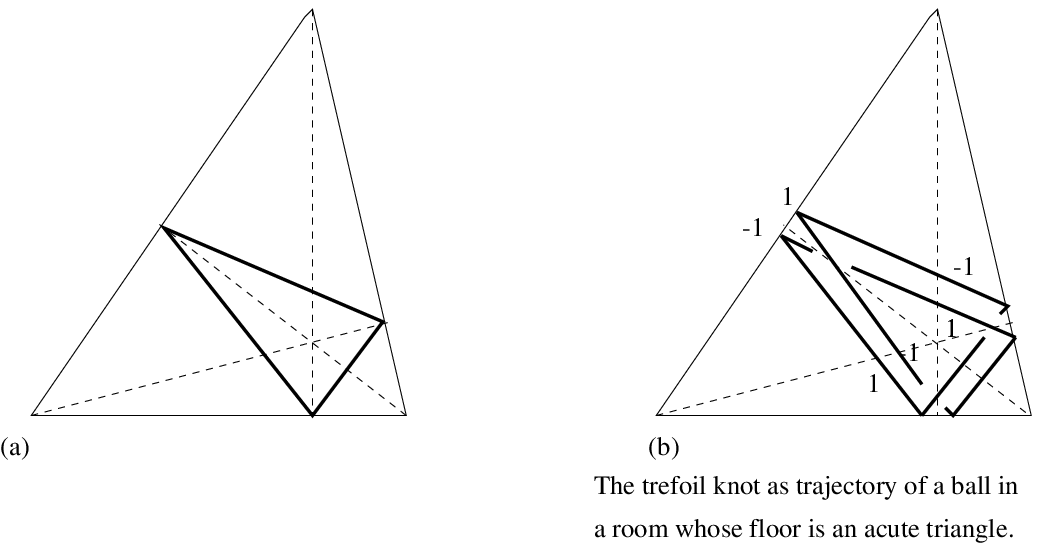}} 
\begin{center} 
Fig. 6.2 
\end{center} 

\item [(ii)] The trivial knot is a trajectory of a ball in a room with an 
right triangular floor, Fig.6.3. 
\ \\ 
\centerline{\psfig{figure=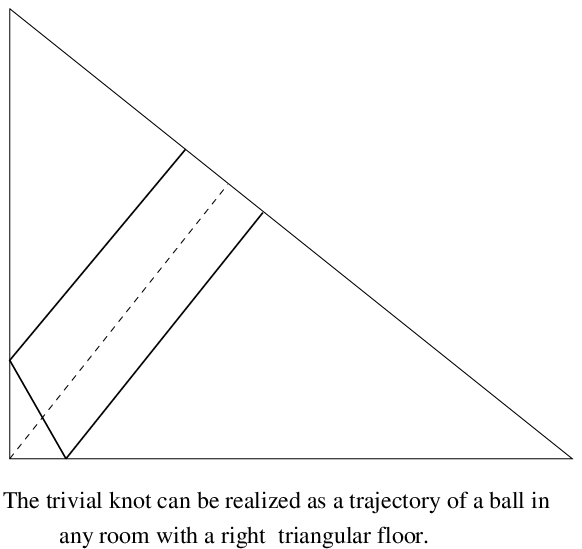}} 
\begin{center} 
Fig. 6.3 
\end{center} 

\item [(iii)] If the floor of a room is a general obtuse triangle, it is an 
open problem whether any knot can be realized as the trajectory of a ball 
in it. However we have general theorem that periodic points are dense 
(in the phase space of the billiard flow) in a rational polygon 
(that is, all polygonal angles are rational with respect to $\pi$) 
 \cite{BGKT}. 
\end{enumerate}
\end{example} 
Example 6.2(i) is of interest because it was shown in \cite{BHJS} 
that the trefoil knot is not a Lissajous knot and thus it is not 
a trajectory of a ball in a room with a rectangular floor. 
More generally we show in Section 3 that no nontrivial torus knot is a 
Lissajous knot. However, we can construct infinitely many 
torus knots in prisms and in the cylinder. 

\begin{example}\label{6.3} 
\begin{enumerate} 
\item [(i)] Any torus knot (or link) of type ($n,2$) can be realized as a 
trajectory of a ball in a room whose floor is a regular $n$-gon ($n\geq 3$). 
Fig.6.1 shows the $(3,2)$ torus knot (trefoil) in the regular triangular 
prism; Fig.6.4(a) depicts the $(4,2)$ torus link in the cube; and 
Fig.6.4(b)(c) illustrates the $(5,2)$ torus knot in a room with 
a regular pentagonal floor. 

\centerline{\psfig{figure=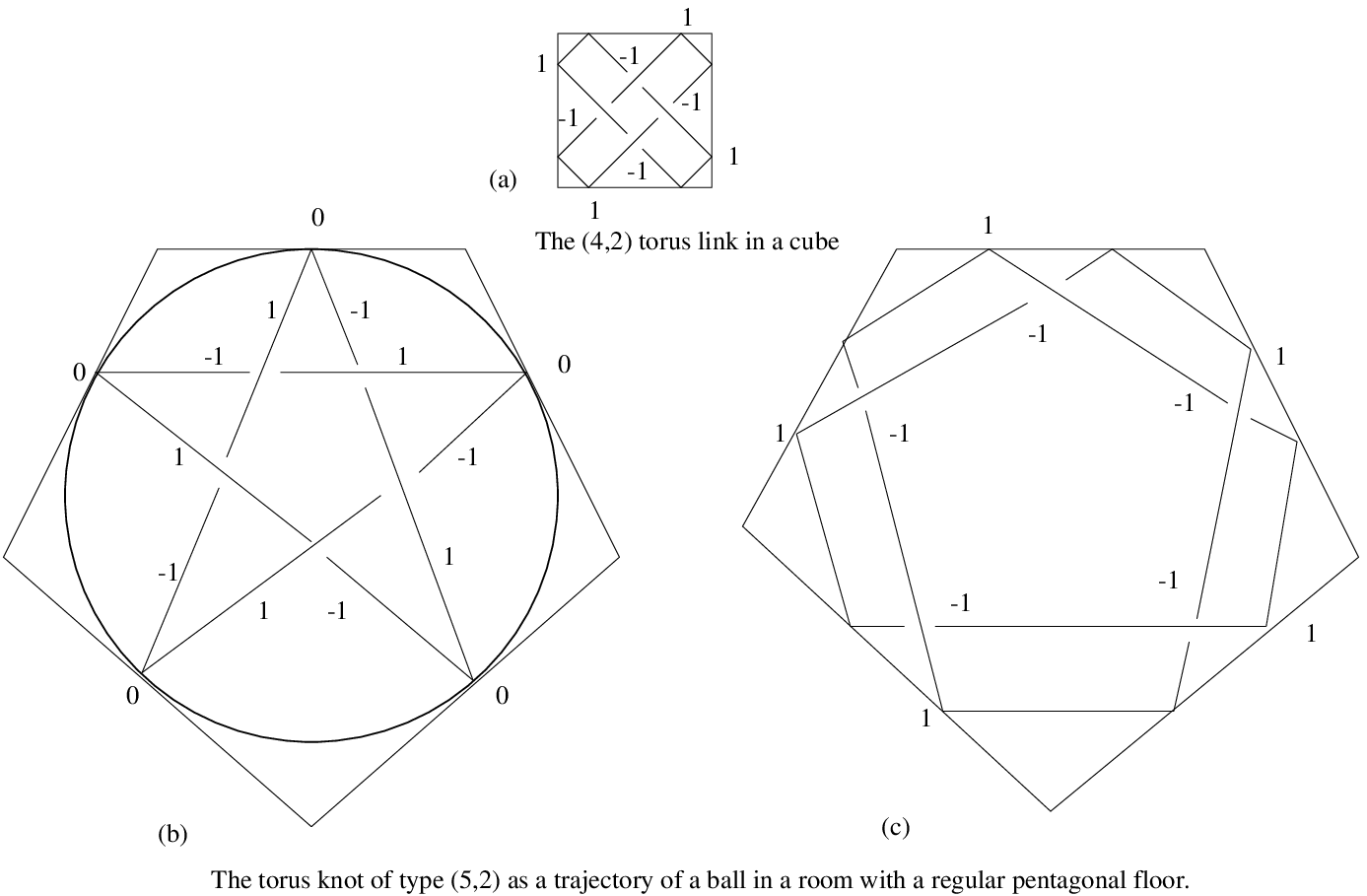}} 
\begin{center} 
Fig. 6.4 
\end{center} 
\item [(ii)] The $(4,3)$ torus knot is a trajectory of a ball in a room with 
the regular octagonal floor; Fig.6.5(a). 

\ \\ 
\centerline{\psfig{figure=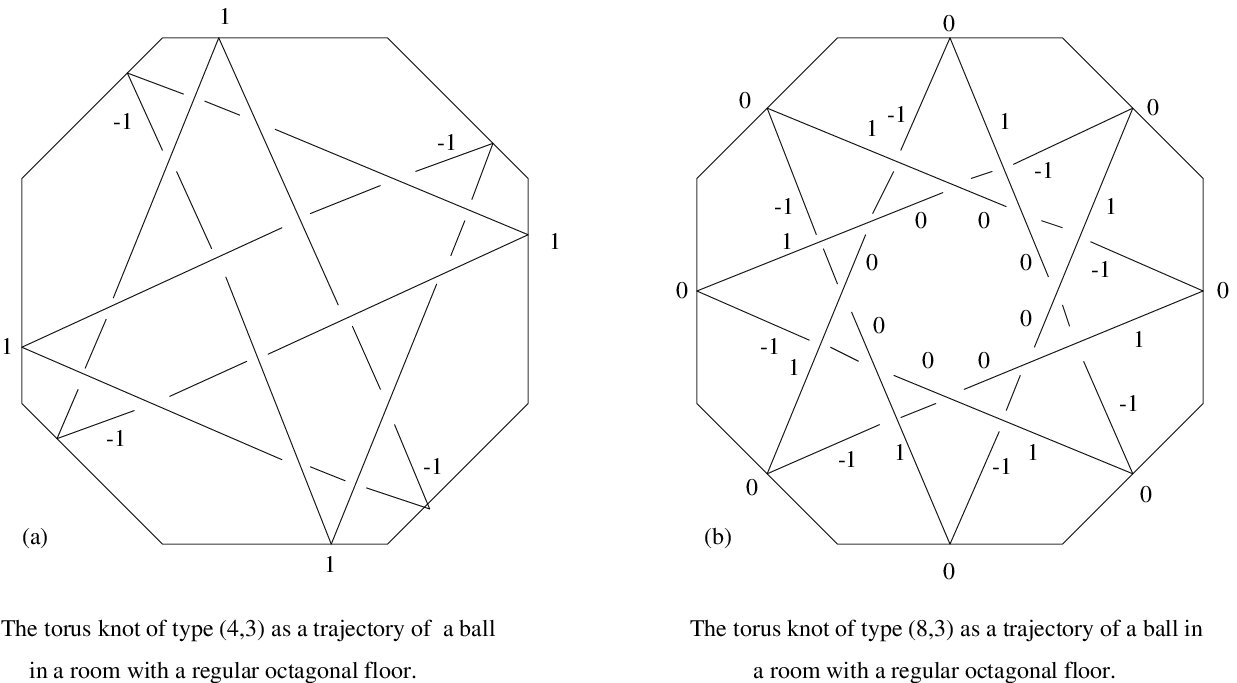}} 
\begin{center} 
Fig. 6.5 
\end{center} 
\ \\ 

\item [(iii)] 
Figures 6.5(b) and 6.6 illustrate how to construct a torus knot (or link) 
of type ($n,3$) in a room with a regular $n$-gonal floor for $n\geq 7$. 
\item [(iv)] 
Any torus knot (or link) of type ($n,k$), where $n\geq 2k+1$, can be 
realized as a trajectory of a ball in a room with a regular $n$-gonal floor. 
The pattern 
generalizes that of Figures 6.4(b), 6.5(b) and 6.6. Edges of the diagram 
go from the center of the $i^{th}$ edge to the center of the $(i+k)^{th}$ edge 
of the n-gon. The ball  bounces from walls at altitude $0$ and its trajectory 
has $n$ maxima and $n$ minima. The whole knot (or link) is $n$-periodic. 
\end{enumerate} 
\end{example} 

\ \\ 
\centerline{\psfig{figure=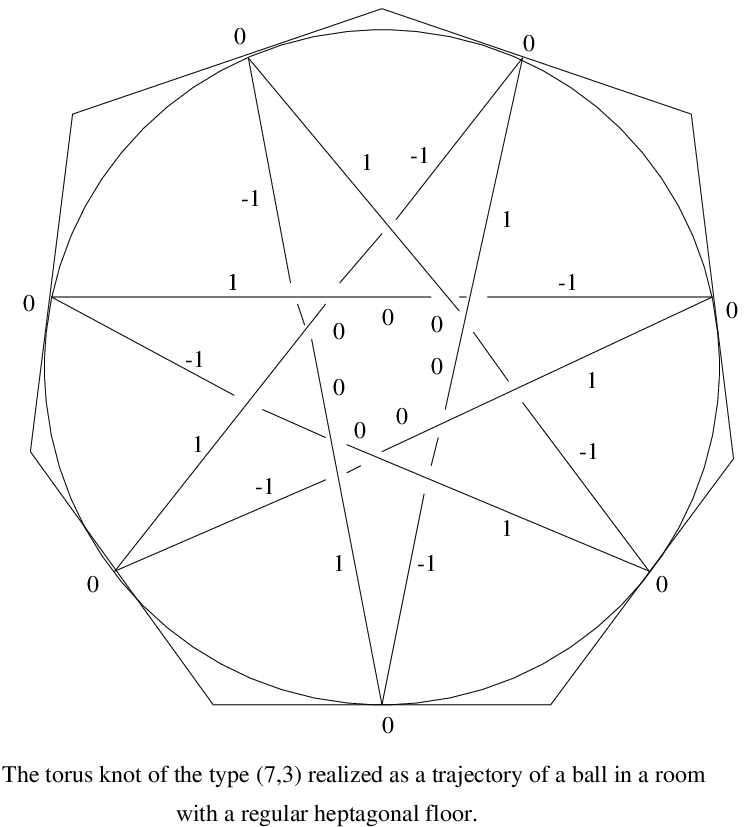}} 
\begin{center} 
Fig. 6.6 
\end{center} 
\ \\ 

\begin{example}\label{6.4} 
Let $D$ be a closed billiard trajectory on a 2-dimensional polygonal table. 
If $D$ is composed of an odd number of segments, then we can always find 
the ``double cover" closed trajectory $D^{(2)}$ in the neighborhood of 
$D$ (each segment will be replaced by two  parallel segments on the 
opposite sides of the initial segment). 
This idea can be used to construct, for a given billiard knot $K$ 
in a polygonal prism (the projection $D$ of $K$ having an odd number 
of segments), a 2-cable $K^{(2)}$ of $K$ as a billiard trajectory 
(with projection $D^{(2)}$). This idea is illustrated in Fig.6.1 and 6.4(c) 
(the (5,2) torus knot as a 2-cable of a trivial one). Starting from 
Example 2.3(iv) we can construct a 2-cable of a torus knot of the 
type ($n,k$) in a regular $n$-gonal prism, for $n$ odd and $n\geq 2k+1$. 
\end{example} 

It follows from \cite{BHJS} that 3-braid alternating knots of 
the form $(\sigma_1\sigma_2^{-1})^{2k}$ are not Lissajous knots as they 
have a non-zero 
\pagebreak
Arf invariant (Corollary 6.12). For $k=1$ we have the 
figure eight knot and for $k=2$ the $8_{18}$ knot \cite{Ro}. 

\begin{example}\label{6.5} 
\begin{enumerate} 
\item [(i)] 
The Listing knot (figure eight knot) can be realized 
as a trajectory of a ball in a room with a regular octagonal floor, 
Fig. 6.7.
\item [(ii)] 
Fig.6.8 describes the knot $8_{18}$ as a trajectory of 
a ball in a room with a regular octagonal floor. This pattern can 
be extended to obtain the knot (or link) which is the closure of the 
three braid $(\sigma_1\sigma_2^{-1})^{2k}$ in a regular $4k$-gonal 
prism ($k>1$). 
\end{enumerate}
\end{example} 
\ \\ 
\centerline{\psfig{figure=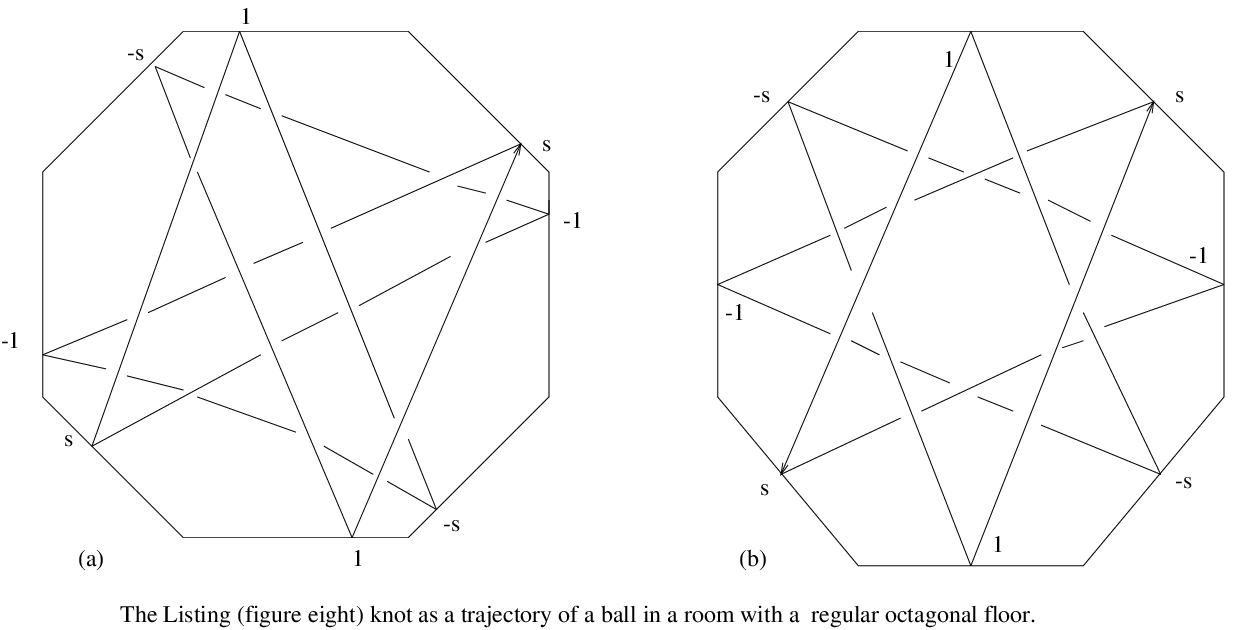}} 
\begin{center} 
Fig. 6.7\footnote{Added for e-print: the journal version of the paper 
illustrates Fig.6.7 by computer graphics generated by Mike Veve.} 
\end{center} 
\ \\ 

\ \\ 
\centerline{\psfig{figure=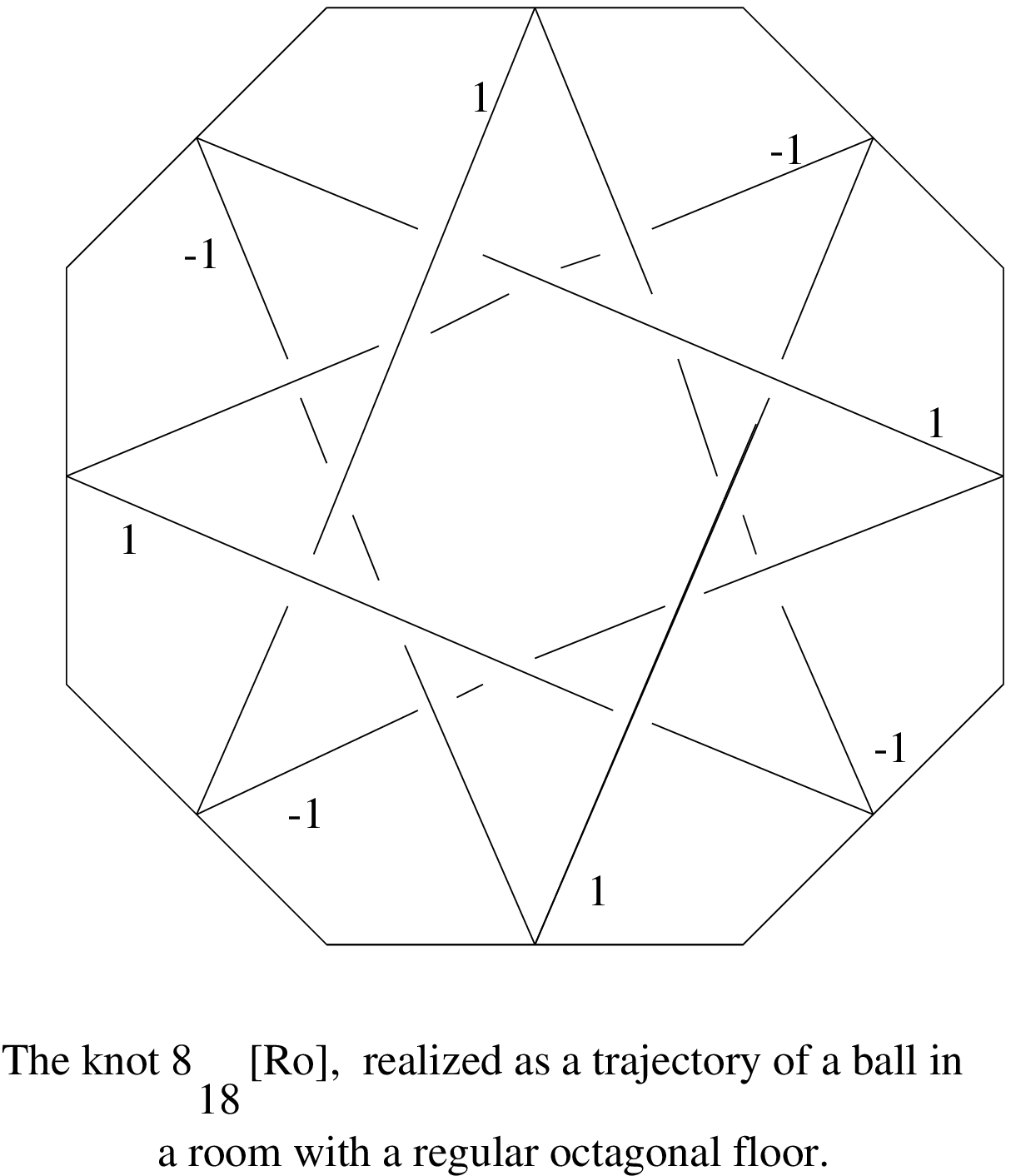,height=6.5cm}} 
\begin{center} 
Fig. 6.8 
\end{center} 
\ \\ 

In the example below we show that the cylinder $D^2 \times [-1,1]$
 support an infinite number of different knot types (in the case of a cube it 
was shown in \cite{La}).

\begin{example}\label{6.6} 
\begin{enumerate} 
\item [(i)] 
Any torus knot (or link) of type ($n,k$), where $n\geq 2k+1$, can be 
realized as a trajectory of a ball in the cylinder; compare Fig.6.4(b), 
Fig.6.5(b) and Fig.6.6. 
\item [(ii)] Every knot (or link) which is the closure of the 
three braid $(\sigma_1\sigma_2^{-1})^{2k}$  can be 
realized as the trajectory of a ball in the cylinder. See Fig. 6.7(b) 
for the case of $k=1$ (Listing knot) and Fig. 6.8 for the case of 
$k=2$ and the general pattern. 
\end{enumerate} 
\end{example} 

Any type of knot can be obtained as 
a trajectory of a ball in some polyhedral billiard (possibly very 
complicated). To see this, consider a polygonal knot in $R^3$ and place 
``mirrors" (walls) at any vertex, in such a way that the polygon is a 
 ``light ray" (ball) trajectory. 
\begin{conjecture}\label{6.7}\ \\ 
Any knot type can be realized 
as the trajectory of a ball in a polytope. 
\end{conjecture} 
\begin{conjecture}\label{6.8}\ \\ 
Any polytope supports an infinite number of different knot types. 
\end{conjecture} 

\begin{problem}\label{6.9}\ \\ 
\begin{enumerate} 
\item [1.] Is there a convex polyhedral billiard in which any knot type 
can be realized as the trajectory of a ball? 
\item [2.] Can any knot type be realized as the trajectory of a ball in a room 
with a regular polygonal floor? 
\item [3.] 
Which knot types can be realized as trajectories of a ball in a 
cylinder ($D^2 \times [-1,1]$)? 
\end{enumerate} 
\end{problem} 
The partial answer to 6.9(3.) was given in \cite{J-P} and \cite{L-O}.
In particular Lamm and Obermeyer have shown that not all knots are
knots in a cylinder (e.g. $5_2$ and $8_{10}$ are not cylinder knots).
The new interesting feature of \cite{L-O} is the use of ribbon condition.

Below we list some information on Lissajous knots (or equivalently
billiard knots in a cube).
\begin{theorem}[\cite{BHJS}]\label{6.10} 
An even Lissajous knot is $2$-periodic and an odd Lissajous knot is 
strongly $+$ amphicheiral. A Lissajous knot is called odd if all, $\eta_x ,
\eta_y$ and $ \eta_z$ are odd. Otherwise it is called an even Lissajous knot.
\end{theorem}

\begin{theorem}[\cite{J-P}]\label{6.11} 
In the even case the linking number of the axis of the $Z_2$-action with the 
knot is equal to $\pm 1$. 
\end{theorem} 
\begin{corollary}\label{6.12} 
\begin{enumerate} 
\item[(i)] (\cite{BHJS}.) 
The Arf invariant of the Lissajous knot is 0. 
\item[(ii)](\cite{J-P}) 
A nontrivial torus knot is not a Lissajous knot. 
\item[(iii)](\cite{J-P}) 
For $\eta_z =2$ a  Lissajous knot is a two bridge knot and its 
Alexander polynomial is congruent to 1 modulo 2. 
\end{enumerate} 
\end{corollary} 
\begin{theorem}[\cite{La}]\label{6.13}
Let $\eta_x, \eta_y > 1$ be relatively prime integers and $K$ the Lissajous
knot with $\eta_z= 2\eta_x \eta_y - \eta_x - \eta_y$ and $\phi_x=
\frac{2\eta_x -1}{2\eta_z}\pi$, $\phi_y= \frac{\pi}{2\eta_z}$, $\phi_z=0$.\\
Then the Lissajous diagram of the projection on the $x-y$ plane is alternating.
Above knots form an infinite family of different Lissajous knots.
\end{theorem}
Motivated by the case of $2$-periodic knots we propose 
\begin{conjecture}\label{6.14}\ \\ 
Turks head knots, (e.g. the closure of the 3-string braids 
$(\sigma_1 \sigma^{-1}_2)^{2k+1}$), are not Lissajous. Observe that they are 
strongly $+$ amphicheiral. 
\end{conjecture} 

We do not think, as the above conjecture shows, that the converse to 
Theorem 6.11 holds. However for $2$-periodic knots it may hold 
(the method sketched in Section 0.4 of \cite{BHJS} may work). 

\begin{problem}\label{6.15}\ \\ 
Let $K$ be a  $Z_2$-periodic knot, such that the linking number of the 
axis of the $Z_2$-action with $K$ is equal to $\pm 1$. 
Is $K$ an even Lissajous knot? 
\end{problem} 
   The first prime knots (in the knot tables \cite{Ro}) which 
may or may not be Lissajous are  $8_3$, $8_6$ ($7_5$ is constructed
in \cite{La}). 

\section{Applications and Speculations}
\ \\ 
A lot can be said about the  structure of a manifold by studying its 
symmetries. 
The existence of $Z_r$ action on a homology sphere is reflected in the 
Reshetikhin-Turaev-Witten  invariants. In our description we follow \cite{K-P}. 
 Let $M$ be a closed connected oriented 3-manifold represented as a 
surgery on a framed link $L\subset S^3$. Let $r\geq 3$, and set the 
variable $A$ used in the Kauffman skein relation to be a primitive 
root of unity of order $2r$. In particular, $A^{2r}=1$. 
Recall  that the invariant ${\mathcal I}_r(M)$ is given by: 

\begin{equation} 
{\mathcal I}_r(M)=\kappa^{-3\sigma_L}\eta^{com(L)}[ L(\Omega_r)] 
\end{equation} 

In $L(\Omega )$ each component of $L$ is decorated by an element $\Omega$ from 
the Kauffman bracket skein module of a solid torus (see Proposition 7.3): 

$$\Omega_r=\sum_{i=0}^{[(r-3)/2]} [e_i] e_i$$ 

Elements $e_i$ satisfy the recursive relation: 
$$e_{i+1}=ze_i-e_{i-1}$$ 
where $z$ can be represented by a longitude of the torus, and $e_0=1$, 
$e_1=z$. 
The value of the Kauffman bracket skein module of the skein element 
$e_i$, when the solid torus is embedded in $S^3$ in a standard way, is 
given by 
$[e_i] = (-1)^i\frac{A^{2i+2}- A^{-2i-2}}{A^2-A^{-2}}$. 
$[e_i]$ is the version of the Kauffman bracket, normalized in such a way
that $[\emptyset]=1$ and $[L]=(-A^2-A^{-2})\langle L \rangle$.
In the equation (2), $\eta$ is a number which satisfies 
$\eta^2[\Omega_r] =1$,  
and $\kappa^3$ is a root of unity such that 
$$\kappa^6 = A^{-6 -r(r+1)/2}.$$ 
Finally, $\sigma_L$ denotes the signature of the linking matrix of 
$L$. 
\begin{theorem}[K-P]\label{7.1}\ \\ 
Suppose that $M$ is a homology sphere and $r$ is an odd prime. If $M$ 
is $r$-periodic then 
$${\mathcal I}_r(M)(A) = 
\kappa^{6j}\cdot {\mathcal I}_r(M)(A^{-1})\mod (r)$$ 
 for some integer $j$. 
\end{theorem} 
Theorem 7.1 holds also for $Z_r$-homology spheres.

Skein modules can be thought as generalizations to 3-manifolds
of polynomial invariants of links in $S^3$. Our periodicity criteria
(especially when link and its quotient under a group action are 
compared), can be thought as the first step toward understanding relations
between skein modules of the base and covering space in a covering.
We will consider here the Kauffman bracket skein module and the relation
between a base and covering space in the case of the solid torus.

\begin{definition} [\cite{P-2,H-P-2}]\label{7.2}\ \\ 
Let $M$ be an oriented 3-manifold, ${\cal L}_{fr}$ the set of unoriented 
framed links in $M$ ($\emptyset$ allowed), 
$R=Z[A^{\pm 1}$, and $R{\cal L}_{fr}$ the free $R$
module with basis ${\cal L}_{fr}$. Let $S_{2,\infty}$ be 
the submodule of $R{\cal L}_{fr}$ generated by skein expressions 
$L_+- AL_0 - A^{-1}L_{\infty}$, where the triple 
$L+, L_0, L_{\infty}$ is presented in Fig.7.1, and $L \sqcup T_1 +(A^2 + 
A^{-2})L$, where $T_1$ denotes the trivial framed knot. 
We define the Kauffman bracket skein module (KBSM), ${\cal S}_{2,\infty}(M)$, 
as the quotient ${\cal S}_{2,\infty}(M)= R{\cal L}_{fr}/S_{2,\infty}$. 
\end{definition} 

\ \\ 
\centerline{\psfig{figure=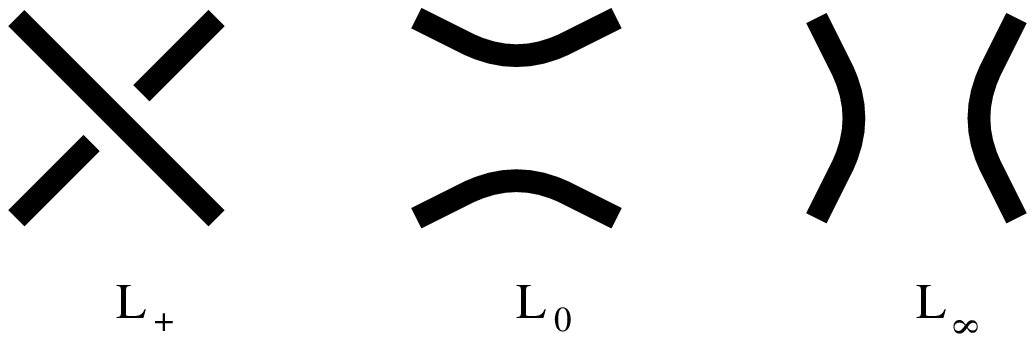,height=3.2cm}} 
\begin{center} 
Fig. 7.1. 
\end{center} 
Notice that $L^{(1)}=-A^3L$ in ${\cal S}_{2,\infty}(M)$; we call this the 
framing relation. In fact this relation can be used instead of 
$L \sqcup T_1 +(A^2 + A^{-2})L$ relation.

\begin{proposition}[\cite{H-P-1,P-2}]\label{7.3}\ \\
The KBSM of a solid torus (presented as an annulus times an interval),
is an algebra generated by a longitude of the solid torus; it is $Z[A^{\pm 1}]$ 
algebra isomorphic to $Z[A^{\pm 1}][z]$ where $z$ corresponds to the longitude.
\end{proposition}

Let the group $Z_r$ acts on the solid torus $(S^1\times [1/2,1])\times [0,1]$) 
with the generator $\varphi (z,t)= (e^{2\pi i/r}z,t)$. 
where $z$ represents an annulus point and $t$ an interval point.
Let $p$ be the $r$-covering map determined by the action (see Fig. 7.2). 
We have the ``transfer" map $p^{-1}_*$ from the KBSM of the base to 
the KBSM of the covering space (modulo $r$) due to the generalization 
of Lemma 3.2. where $z$ represent an annulus point and $t$ an interval point. 
To present the generalization in the  most natural setting we 
introduce the notion of the ``moduli" equivalence of links and associated
moduli skein modules.\ \\
\ \\   
\centerline{\psfig{figure=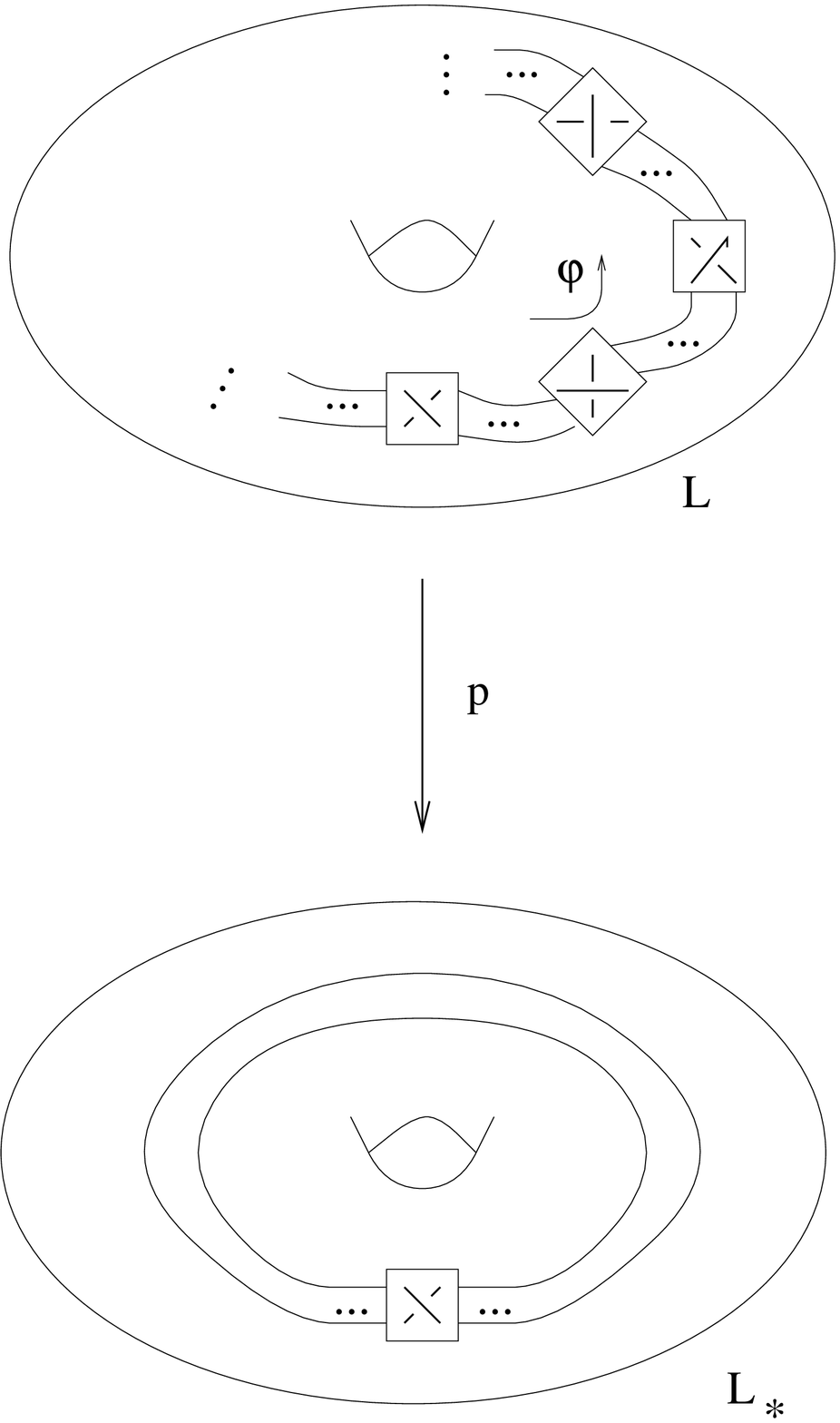,height=8.5cm}} 
\begin{center} 
Fig. 7.2 
\end{center} 
\ \\ 
\begin{definition}\label{7.4}
\begin{enumerate}
\item[(i)] We say that two links in a manifold $M$ are moduli equivalent
if there is a preserving orientation homeomorphism of $M$ sending
one link to another.
\item[(ii)] Let $G$ be a finite group action on $M$. We say that two links 
$L$ and $L'$ in $M$ are $G$ moduli equivalent if $L'$ is ambient isotopic to
$g(L)$ for some $g\in G$.
\item [(ii)] We define moduli KBSM (resp. $G$ moduli KBSM) as KBSM divided
by moduli (resp. $G$ moduli) relation. We would use the notation
${\cal S}^m_{2,\infty}(M)$ (resp. ${\cal S}^G_{2,\infty}(M)$).
\end{enumerate}
\end{definition}

\begin{lemma}\label{7.5} 
Let  the group $Z_r$, $r$ prime, act on the oriented 3-manifold $M$ and let
 $L=L_{sym(cr)}$ be a framed singular link in $M$ which satisfies 
$\varphi (L) = L$  where $\varphi: M \to M$ is the generator of the 
$Z_r$-action, and $Z_r$ has no fixed points on $L$. 
Then in the skein module 
${\cal S}^{Z_r}_{2,\infty}(M)$ one has the formula 
$$ L_{sym ({\psfig{figure=L+nmaly.eps}})} \equiv 
A^r L_{sym ({\psfig{figure=L0nmaly.eps}})} + 
A^{-r} L_{sym ({\psfig{figure=Linftynmaly.eps}})} \mod r$$ 
where 
$L_{sym( 
{\psfig{figure=L+nmaly.eps}}) 
}, L_{sym({\psfig{figure=L0nmaly.eps}} 
)}$ and $L_{sym({\psfig{figure=Linftynmaly.eps}} 
)}$ 
denote three $\varphi$-invariant diagrams of links 
which are the same outside of the $Z_r$-orbit of a neighborhood of a
fixed singular crossing at which they differ by replacing
$L_{sym(cr)}$ by
{\psfig{figure=L+nmaly.eps}} 
or
{\psfig{figure=L0nmaly.eps}}, 
or 
{\psfig{figure=Linftynmaly.eps}}, 
respectively.\\
Notice that in the case of an action on the solid torus, $\varphi$ is isotopic
to identity, thus $Z_r$ moduli KBSM is the same as KBSM. 

\end{lemma} 

We can use Lemma 7.5 to find the formula for a torus knot in a solid torus
 modulo $r$, for a prime $r$. 
\begin{theorem}\label{7.6}  The torus knot $T_{r,k}$ satisfies:
$$T(r,k) = A^{r(k-1)}(\frac{x^{k+1}-x^{-k-1}-A^{-4r}(x^{k-1}-x^{1-k})} 
 {x-x^{-1}})$$ 
where $z=x+x^{-1}$ is a longitude of the solid torus (an annulus times 
an interval). 
\end{theorem}
One can prove Theorem 7.6 by rather involved 
induction\footnote{Added for e-print: The formula has been 
generalized and given a simple
ideological proof in R.Gelca, C.Frohman, 
Skein Modules and the Noncommutative
Torus, {\it Transactions of the AMS}, 352, 2000, 4877-4888.} 
(compare example 7.7),
however modulo $r$ Theorem 7.6 has very easy proof using Lemma 7.5.
Namely $T_{r,k}$ is an $r$ cover of a knot $T_{1,k}$ and for $T_{1,k}$
one has $T_{1,k+2}= A(x+x^{-1})T_{1,k+1} - A^2T_{1,k}$ in KBSM, thus
$T_{1,k}= A^{(k-1)}(\frac{x^{k+1}-x^{-k-1}-A^{-4}(x^{k-1}-x^{1-k})} 
 {x-x^{-1}})$
Now one compares binary computational trees of $T_{1,k}$ and $T_{r,k}$ 
(modulo $r$) using lemma 7.5, to get $\mod r$ version of Theorem 7.6.\\
In fact for any regular $r$ covering $p:M \to M_*$ we have the transfer
map (\mod $r$):
$$p^{-1}_*: {\cal S}_{2,\infty}(M_*)  \to {\cal S}^{Z_r}_{2,\infty}(M) \mod r$$
which is a $Z$ homomorphism and $p^{-1}_*(w(A)L) \equiv w(A^r)p^{-1}(L) \mod r$.
\begin{example}\label{7.7} 
We compute here the torus link $T_{2,k}$ in the KBSM of a solid torus. 
We will prove that for an odd $k$ (that is when $T_{2,k}=K_{2,k}$ is a knot): 
$$ L(2,k) = A^{2(k-1)}\frac{x^{k+1}-x^{-k-1}- 
A^{-8}(x^{k-1}- x^{1-k})}{x-x^{-1}}$$ 
and for an even $k$: 
$$ L(2,k) = A^{2(k-1)}\frac{x^{k+1}-x^{-k-1}- 
A^{-8}(x^{k-1}- x^{1-k})}{x-x^{-1}} + 2A^{-6}$$ 
The formulas follow inductively from the following, easy to check, 
recurrence relation ($k>2$), compare Fig 7.3.: 
$$L(2,k) = A^2L_{-,-}^{a,b} - A^{-2}L_{|,|}^{a,b} + A^{-1}L_{+,|}^{a,b} + 
A^{-1}L_{|,+}^{a,b}= A^2zL_{2,k-1} - A^4L_{2,k-2} -2A{-4}u_k$$ where 
$u_k= z$ for $k$ odd and $u_k=-A^2 - A^{-2}$ for $k$ even. 
If we denote $e_n=\frac{x^{n+1}-x^{-n-1}}{x-x^{-1}} $, Chebyshev polynomial 
in the variable $z=x+x^{-1}$, then in odd cases the formula can be written as: 
$$ L(2,k) = A^{2(k-1)}(e_k- A^{-8}e_{k-2})$$. 
\end{example} 

\ \\ 
\centerline{\psfig{figure=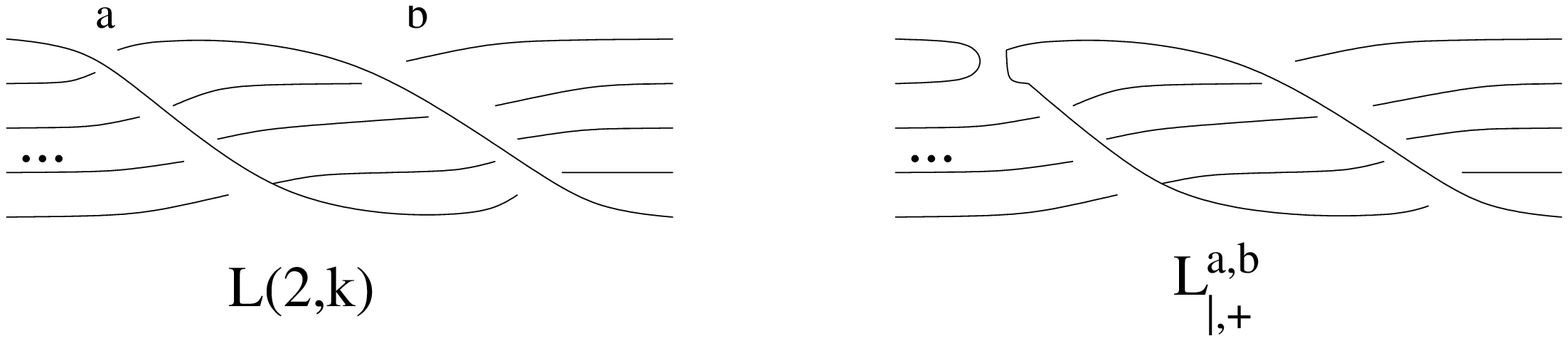,height=2.6cm}} 
\begin{center} 
Fig. 7.3
\end{center} 

\begin{corollary}[Traczyk]\label{7.8} \ \\
\begin{enumerate} 
\item [(a)] $A^3e_3$ is represented in the solid torus 
by a 2-component link being a closure of the 3-braid $\Delta = 
\sigma_1 \sigma_2 \sigma_1$; Fig. 7.4. 
\item [(b)] In $S^3$, the closed braid links $\hat{\Delta}^{2n+1}$ form an 
infinite family of (simple) links whose Jones polynomial differ only 
by an invertible element\footnote{Added for e-print:
P.Traczyk, 3-braids with proportional Jones polynomials,
{\it Kobe J. Math.}, 15(2), 1998, 187-190.}
 in $Z[t^{\pm \frac{1}{2}}]$. 
\item [(c)] There are infinite families of links sharing the same 
Jones polynomial (e.g. $t^2 + 2 + t^{-2}$)). 
\end{enumerate} 
\end{corollary} 
\ \\ 
\centerline{\psfig{figure=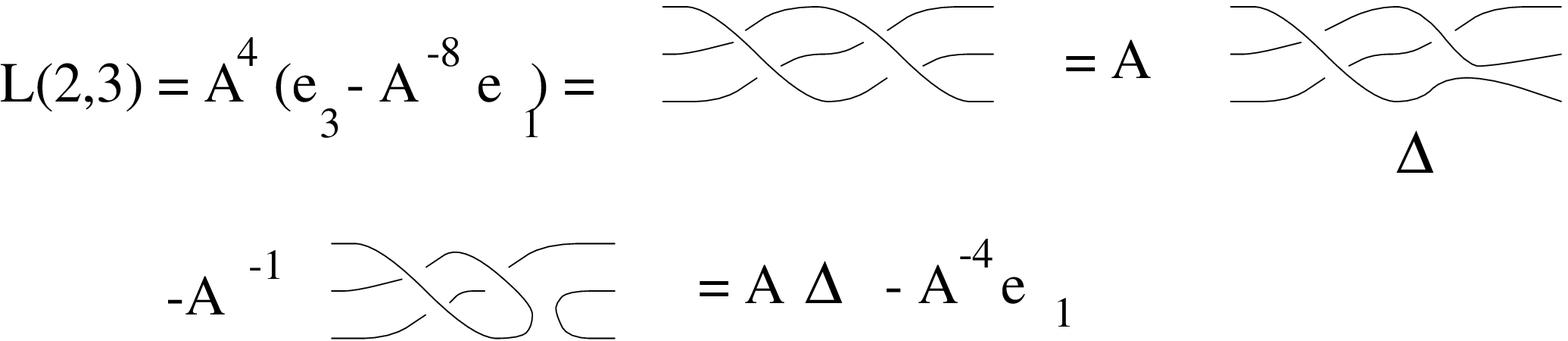,height=2.6cm}} 
\begin{center} 
Fig. 7.4
\end{center} 

\begin{proof} 
\begin{enumerate} 
\item [(a)] Using the formula of Example we have 
$L_{2,3}= A^4(e_3 - A^{-8}e_1)$. The necessary calculation is shown 
in Fig. 7.4. 
\item [(b)] 
$e_3$ is an eigenvector of the Dehn twist on the solid torus and 
 any braid, $\gamma$, is changed by a Dehn twist to ${\Delta}^2 \gamma$. 
Thus $\hat{\Delta}^{2n+1}= A^{6n}\hat{\Delta}$ in the KBSM of the solid torus.\\ 
Closed braids $\hat{\Delta}^{2n+1}$ form an infinite family of different 
2-component links in $S^3$ (linking number equal to $k$). 
\item [(c)] Consider a connected sum of $\hat{\Delta}^{2n+1}$ and 
$\hat{\Delta}^{2m+1}$. For fixed $n+m$ one has an infinite family 
of different 3-component links with the same Jones polynomial. 
For $m=-n-1$ one get the family of links with the Jones polynomial of the 
connected sum of right handed and left handed Hopf links (thus 
$(t+t^{-1})^2$). 
\end{enumerate} 
\end{proof} 
We do not know which $e_n$ can be realized by framed links in a solid torus,
however we have:
\begin{lemma}\label{7.9}
If $e_n$ is realized by a framed link then $n=2^k-1$ for some $k$.
\end{lemma} 
\begin{proof}
Consider the standard embedding of a solid torus in $S^3$. Then
$<e_n>_{A=-1} = (-1)^{n-1}\frac{n+1}{2}$ on the other hand for a link
$L$ one has $<L>_{-1} = (-2)^{com(L)-1}$. Thus $n=2^k-1$
\end{proof}

Address:\\ 
Department of Mathematics\\ 
George Washington University  \\ 
Washington, DC 20052 \\ 
e-mail: przytyck@math.gwu.edu\\
\{przytyck@gwu.edu\} 

\end{document}